\documentclass[twoside,11pt]{article}

\usepackage{amssymb, amsmath, latexsym, mathrsfs, verbatim, calc}

\def \N{{\rm I\!N}}

\newcommand{\dem}{\noindent \underline{Proof} :}
\newcommand{\cq}{\hfill$\Box$\\}

\def\i{{\bf 1}}

\def \R{{\mathop{{\rm I\negthinspace R}}}}

\def\to{\rightarrow}

\def\bart{{\bar t}}
\def\bars{{\bar s}}
\def\barx{{\bar x}}
\def\bary{{\bar y}}
\def\barp{{\bar p}}
\def\barq{{\bar q}}
\def\hp{{\hat p}}
\def\hq{{\hat q}}
\def\hw{{\hat w}}

\newcommand{\ba}{\begin{array}}
\newcommand{\ea}{\end{array}}
\newcommand{\be}{\begin{equation}}
\newcommand{\ee}{\end{equation}}
\newcommand{\bea}{\begin{eqnarray}}
\newcommand{\eea}{\end{eqnarray}}
\newcommand{\beaa}{\begin{eqnarray*}}
\newcommand{\eeaa}{\end{eqnarray*}}
\newtheorem{thm}{Theorem}[section]
\newtheorem{Lemma}[thm]{Lemma}

\newtheorem{prop}[thm]{Proposition}

\newtheorem{Proposition}{Proposition}[section]

\newtheorem{Corollary}{Corollary}[section]

\newtheorem{Definition}{Definition}[section]
\newcommand{\AR}{{\cal A}}
\newcommand{\BR}{{\cal B}}
\newcommand{\CR}{{\cal C}}

\newcommand{\FR}{{\cal F}}

\newcommand{\PR}{{\cal P}}
\newcommand{\UR}{{\cal U}}
\newcommand{\VR}{{\cal V}}
\newcommand{\vme}{V^{-*}}
\newcommand{\vome}{V^{-*}_1}
\newcommand{\ar}{\AR_r}
\newcommand{\br}{\BR_r}

\begin{document}

\title{Stochastic differential games with asymmetric information.}

\author{Pierre Cardaliaguet and Catherine Rainer\\
$\;$\\
\small D\'epartement de Math\'ematiques, Universit\'e de Bretagne Occidentale,
\\\small 6, avenue Victor-le-Gorgeu, B.P. 809, 29285 Brest cedex, France\\
e-mail : 
Pierre.Cardaliaguet@univ-brest.fr, Catherine.Rainer@univ-brest.fr}
  \maketitle

\noindent {\bf Abstract : }  We investigate a two-player zero-sum stochastic differential game in which the players have an
asymmetric information on the random payoff. We prove that the game has a value and characterize this 
value in terms of {\it dual} solutions of some second order Hamilton-Jacobi equation.
\vspace{3mm}

\noindent{\bf Key-words : } stochastic differential game, asymmetric information, viscosity solution.
\vspace{3mm}

\noindent{\bf A.M.S. classification :} 49N70, 49L25, 91A23.
\vspace{3mm}

\section{Introduction}

This paper is devoted to a class of two-player zero-sum stochastic differential game in which the players have different information
on the payoff.  In this basic model, the terminal cost is chosen (at the initial time) randomly among a finite set of costs  $\{g_{ij}, \; i\in \{1, \dots, I\}, \; j\in \{1, \dots, J\}\,\}$. 
More precisely, the indexes $i$ and $j$ are chosen independently according to a probability $p\otimes q$ on 
$\{1, \dots, I\}\times\{1, \dots, J\}$. 
Then the index $i$ is announced to the first player and  the  index $j$ to the second player.
The players control the stochastic differential equation
$$
\begin{array}{l}
dX_s=b(s,X_s,u_s,v_s)ds+\sigma(s,X_s,u_s,v_s)dB_s,\; s\in[t,T],\\
X_t=x,
\end{array}
$$
through their respective controls $(u_s)$ and $(v_s)$ in order, for the first player, to minimize $E[g_{ij}(X_T)]$ and, for 
the second player, to maximize this quantity. Note that the players do not really know which payoff they are
actually  optimizing because the first player, for instance, ignores which index $j$ has been chosen. The key assumption in our
model is that the players observe the evolving state $(X_s)$.
So they can deduce from this observation the behavior of their opponent and try to derive from it some knowledge on their missing data. 

The formalization of such a game is quite involved: we refer to the second section of the paper where the notations are properly
defined. In order to describe our results, let us introduce the upper and lower value functions $V^+$ and $V^-$ of the game:
$$
V^+(t,x,p,q)=\inf_{\hat\alpha\in(\AR_r(t))^I}\sup_{\hat\beta\in
\BR_r(t))^J}J^{p,q}(t,x,\hat\alpha,\hat\beta),
$$
$$
V^-(t,x,p,q)=\sup_{\hat\beta\in
(\BR_r(t))^J}\inf_{\hat\alpha\in(\AR_r(t))^I}J^{p,q}(t,x,\hat\alpha,\hat\beta).
$$
where $J^{p,q}(t,x,\hat\alpha,\hat\beta)$ is the expectation under the probability $p\otimes q$ of
the payoff associated with the strategies $\hat \alpha=(\alpha_i)_{i\in \{1, \dots, I\}}$ and $\hat \beta=(\beta_j)_{j\in \{1, \dots, J\}}$ of the players. 
The strategy $\hat \alpha$ takes into account the
knowledge by the first player of the index $i$  while $\hat \beta$ takes into account the knowledge of $j$ by the second player. 
Our main result is that, under Isaacs'condition, the two value functions coincide: $V^+=V^-$.
Moreover, ${\bf V}:=V^+=V^-$ is the unique viscosity solution {\it in the dual sense}
of some second order Hamilton-Jacobi equation. This means that 
\begin{itemize}
\item[(i)] ${\bf V}$ is convex with respect to $p$ and concave with respect to $q$,
\item[(ii)] the convex conjugate of ${\bf V}$ with respect to $p$ is a subsolution of some Hamilton-Jacobi-Isaacs (HJI) equation in the viscosity sense,
\item[(iii)] the concave conjugate of ${\bf V}$ with respect to $q$ is a supersolution of a symmetric HJI equation,
\item[(iv)] ${\bf V}(T,x,p,q)=\sum_{i,j} p_iq_j g_{ij}(x)$ where $p=(p_i)_{i\in\{1, \dots, I\}}$ and $q=(q_j)_{j\in\{1, \dots, J\}}$.
\end{itemize}
We strongly underline that in general the value functions are {\it not } solution of the standard HJI equation: indeed ${\bf V}$ does not satisfy a dynamic
programming principle in a classical sense.

An important current in Mathematical Finance is the modeling of insider trading (see for example Amendinger, Becherer, Schweizer \cite{abs} or Corcuera, Imkeller, Kohatsu-Higa, Nualart \cite{cikn} and references therein). The basic question studied in these works is to evaluate how the addition of knowledge for a trader---i.e., mathematically, 
the addition to the original filtration of a variable depending on the future---shows up in his investing strategies, and an important tool is the theory of enlargement of filtrations.
Our approach is completely different.  Indeed, what is important in our game 
is not that the players have ``more" information than what is contained in the filtration of the Brownian motion, but that their information differs from that of their opponent.  
In some sense we try to understand the strategic role of information in the game.

The model described above is strongly inspired by a similar one studied by Aumann and Maschler in the framework of repeated games.
Since their seminal papers (reproduced in  \cite{AM}), this model has attracted a lot of attention in game theory (see 
\cite{DM1}, \cite{MZ}, \cite{RSV}, \cite{So}). However
it is only recently that the first author has adapted the model to deterministic differential games  (see \cite{c1}, \cite{c2}). 

The aim of this paper is to generalize the results of \cite{c1}  to stochastic differential games and to game with integral 
payoffs. There are several difficulties towards this aim. First the notion of strategies for stochastic differential
games is quite intricated (see \cite{FS}, \cite{NI}). For our game it is all the more difficult that the players have 
to introduce additional noise in their strategies in order to  confuse their oponent. One of the achievements of this paper is an
important simplification of the notion of strategy which allows the introduction of the notion of random strategies. 
This also simplifies several proofs of \cite{c1}. 
Second the existence of a value for ``classical" stochastic differential games relies on a comparison principle for some second order Hamilton-Jacobi equations.
Here we have to be able to compare functions satisfying the condition (i,ii,iv) defined above 
with functions satisfying (i,iii,iv). While for deterministic differential games (i.e., first order HJI equations)
we could do this without too much trouble (see \cite{c1}), for stochastic differential games (i.e., second order HJI equations) the proof is much more involved. In particular
it requires a new maximum principle for lower semicontinuous functions (see the appendix) which is  the most technical part of the paper. 

The paper is organized in the following way: in section 2, we introduce the main notations and  the notion of 
random strategies and we define the value functions of our game. In section 3 we prove that the value functions (and more precisely the convex and concave
conjugates) are sub- and supersolutions of some HJ equation. 
Section 4 is devoted to the comparison principle and to the existence of the value. In Section 5 we investigate 
stochastic differential games with a running cost. 
The appendix is devoted to a new maximum principle.

\section{Definitions.}

\subsection{The dynamics.}

Let $T>0$ be a fixed finite time horizon.
For $(t,x)\in[0,T]\times\R^n$, 
we consider the following doubly controlled stochastic
system :
\begin{equation}\label{dyn}
\begin{array}{l}
dX_s=b(s,X_s,u_s,v_s)ds+\sigma(s,X_s,u_s,v_s)dB_s,\; s\in[t,T],\\
X_t=x,
\end{array}
\end{equation}
where 
$B$ is a $d$-dimensional standard Brownian motion on a given probability space 
$(\Omega, \FR,P)$.
For $s\in [t,T]$, we set
\[ \FR_{t,s}=\sigma\{ B_r-B_t, r\in [t,s]\}\vee\PR,\]
where $\PR$ is the set of all null-sets of $P$.

 The processes $u$ and  $v$ are assumed to take
their values in some compact metric spaces $U$
and $V$ respectively. 
We suppose 
that the functions $b:[0,T]\times
\R^n\times U\times V\to\R^n$
and $\sigma:[0,T]\times \R^n\times U\times V\to\R^{n\times d}$ are 
continuous and satisfy the
assumption (H):\\

\indent (H) $b$ and $\sigma$ are bounded and Lipschitz
continuous with respect
to $(t,x)$,
uniformly in $(u,v)\in U\times V$.\\

We also assume  Isaacs' condition : for all  $(t,x)\in[0,T]\times
\R^n$, $ p\in\R^n$,  and all $A\in {\cal S}_n$
(where ${\cal S}_n$ is the set of symmetric $n\times n$ matrices) holds:
\begin{equation}\label{Isaacs}
\begin{array}{l}
\inf_u\sup_v\{ <b(t,x,u,v),p>+\frac12
Tr(A\sigma(t,x,u,v)\sigma^*(t,x,u,v))\}=\\
\hspace{1cm}\sup_v\inf_u\{ <b(t,x,u,v),p>+\frac12
Tr(A\sigma(t,x,u,v)\sigma^*(t,x,u,v))\}
\end{array}
\end{equation}
We set $H(t,x,p,A)=\inf_u\sup_v\{ <b(t,x,u,v),p>+\frac12
Tr(A\sigma(t,x,u,v)\sigma^*(t,x,u,v))\}$.

For $t\in [0,T)$, we denote by $\CR([t,T],\R^n)$ the set of continuous maps from $[t,T]$ to $\R^n$.

\subsection{Admissible controls.}

\begin{Definition} An {\em admissible  control} $u$ for player
I (resp. II) on $[t,T]$
is a process taking
values in $U$
(resp. $V$), progressively measurable with respect to the filtration
$(\FR_{t,s}, s\geq t)$. \\
The set of admissible controls for player I (resp. II) on
$[t,T]$
is denoted by ${\cal U}(t)$ (resp. ${\cal V}(t)$).
\end{Definition}

We identify two processes $u$ and $\overline{u}$ in $\UR(t)$ if
$P\{u=\overline u\mbox{ a.e. in }[t,T]\}=1$.\\

Under  assumption (H), for all $(t,x)\in [0,T]\times\R^n$ and
$(u,v)\in\UR(t)\times\VR(t)$, there exists a unique solution to
(\ref{dyn}) that we denote by $X^{t,x,u,v}_.$.

\subsection{Strategies.}

\begin{Definition} 
\label{strat}
A {\em strategy} for player I starting at time $t$ is a Borel-measurable map 
$\alpha : [t,T]\times\CR([t,T],\R^n)\rightarrow U$ 
for which there exists $\delta>0$ such that,
$\forall s\in[t,T], f,f'\in\CR([t,T],\R^n)$, 
if $f=f'$ on $[t,s]$, then $\alpha(\cdot,f)=\alpha(\cdot,f')$ on  $[t,s+\delta]$.\\
We define strategies for player II in a symmetric way and denote by $\AR(t)$ (resp. 
$\BR(t)$) the set of strategies for player I (resp. player II).
\end{Definition}

We have the following existence result :

\begin{Lemma}\label{ptfix}
For all $(t,x)$ in $[0,T]\times\R^n$, for all $(\alpha,\beta)\in\AR(t)\times\BR(t)$, 
there exists a unique couple of controls $(u,v)\in\UR(t)\times\VR(t)$ that satisfies $P-$a.s.
\begin{equation}
 \label{alphau}
(u,v)=(\alpha(\cdot,X^{t,x,u,v}_\cdot),\beta(\cdot,X^{t,x,u,v}_\cdot)) 
\mbox{ on }
[t,T].
\end{equation}

\end{Lemma}

\noindent{\em Proof:} The controls $u$ and $v$ will be built step by step. Let $\delta>0$ 
be a common delay for $\alpha$ and $\beta$. We can choose $\delta$ such that $T=t+N\delta$ for some $N\in\N^*$.\\
By definition, on $[t,t+\delta)$, for all $f\in\CR([t,T],\R^n)$, $\alpha(s,f)=\alpha(s,f(t))$. 
Since, for all $(u,v)\in\UR(t)\times\VR(t)$, $X^{t,x,u,v}_t=x$, 
the control $u$ is uniquely defined on $[t,t+\delta)$ by
\[ \forall s\in [t,t+\delta),u(s)=\alpha(s,x).\]
The same holds for $v$, what permits us to define the process $X^{t,x,u,v}_\cdot$ on $[t,t+\delta)$ as
a solution of the system (\ref{dyn}) restricted on the interval $[t,t+\delta)$.\\
Now suppose that $u$, $v$ and $X^{t,x,u,v}_\cdot$ are $P-$a.s. defined uniquely on some interval $[t,t+k\delta)$, 
$k\in\{ 1,\ldots,N-1\}$.
This allows us to set,
\[ \forall s\in [t+k\delta,t+(k+1)\delta), \; 
u_s=\alpha(s,X^{t,x,u^k,v^k}_\cdot),
v_s=\beta(s,X^{t,x,u^k,v^k}_\cdot),\]
where
\[(u^k,v^k)=
\left\{\begin{array}{l}
(u,v)\mbox{ on }[t,t+k\delta)\\
(u_0,v_0)\mbox{ else,}
\end{array}\right.
\]
for some arbitrary $(u_0,v_0)\in\UR(t)\times\VR(t)$.\\
Considering $X^{t,x,u^k,v^k}_\cdot$ as a random variable 
with values in the set of paths
$\CR([t, T), \R^n)$, 
it is clear that the map $(s,\omega)\to u_s(\omega)$ (defined on $[t+k\delta,t+(k+1)\delta)\times \Omega$)  
as the composition of the Borel measurable application $\alpha$ with the map $(s,\omega)\to (s,X^{t,x,u^k,v^k}_\cdot(\omega))$, is a process on 
$[t+k\delta,t+(k+1)\delta)$ with measurable paths. 
Further, the non anticipativity of $\alpha$ guaranties that, for all $s\in [t+k\delta,t+(k+1)\delta)$,
 $u_s$ 
is $\FR_{t,t+k\delta}$-measurable and the process $u|_{[t,t+(k+1)\delta)}$ is $(\FR_{t,s})$-progressively measurable.
The same holds of course for $v|_{[t,t+(k+1)\delta)}$.\\
With $(u,v)$ defined on $[t,t+(k+1)\delta)$, we can now 
define the process $X^{t,x,u,v}_\cdot$ up to time $t+(k+1)\delta$. 
This completes the proof by induction.\cq

We denote by $X^{t,x,\alpha,\beta}_\cdot$ the process $X^{t,x,u,v}_\cdot$, with
$(u,v)$ associated to $(\alpha,\beta)$ by relation (\ref{alphau}).\\

In the frame of incomplete information it is necessary to introduce random strategies. 
In contrast with \cite{c1} and \cite{c2}, where the
random probabilities are supposed to be absolutely continuous 
with respect to the Lebesgue measure,
play a random strategy will consist here  to choose some strategy in a finite set of 
possibilities, i.e. the involved probabilities are finite. It is not clear if this assumption
is more realistic nor if the notation will be lighter, nevertheless this alternative
allows us to avoid some technical steps of measure theory, 
in a paper that is already technical enough.\\

\noindent{\em Notation:}
 For $R\in\N^*$, let $\Delta(R)$ be the set of all $(r_1,\ldots,r_R)\in[0,1]^R$ 
that satisfy $\sum_{n=1}^Rr_n=1$.\\
We define a random strategy $\overline\alpha$ for player I by
$\overline\alpha=(\alpha^1,\ldots \alpha ^R;r^1,\ldots,r^R)$, with $R\in\N^*$, 
$(\alpha^1,\ldots \alpha^R)\in(\AR(t))^R,\; (r^1,\ldots ,r^R)\in\Delta(R)$.\\
The heuristic interpretation of  $\bar \alpha$ is that player I's strategy amounts to
choose the pure strategy $\alpha^k$ with probability $r^k$. \\
We define in a similar way the random strategies for player II, and denote by 
$\AR_r(t)$ (resp. $\BR_r(t)$) the set of all random strategies for player I (resp. player II).\\
Finally, identifying $\alpha\in\AR(t)$ with $(\alpha;1)\in\AR_r(t)$, 
we can write $\AR(t)\subset\AR_r(t)$, and the same holds for $\BR(t)$ and $\BR_r(t)$.\\

\subsection{The payoff.}

Fix $I,J\in\N^*$.\\
For $1\leq i\leq I, 1\leq j\leq J$, let $g_{ij}:\R^n\to \R$ be the terminal payoffs.  We assume that
\begin{equation}\label{Hypgij}
\mbox{\rm For $1\leq i\leq I, 1\leq j\leq J$, $g_{ij}$ are Lipschitz continuous and bounded.}
\end{equation}
For $(p,q)\in\Delta(I)\times\Delta(J)$, with $p=(p_1,\ldots, p_I),\;
q=(q_1,\ldots q_J)$, we denote with a hat the elements of $(\AR_r(t))^I$ 
(resp. $(\BR_r(t))^J$):
$\hat\alpha=(\overline\alpha_1,\ldots,\overline\alpha_I)$,
$\hat\beta=(\overline\beta_1,\ldots,\overline\beta_J)$.\\
We adopt following notations :\\
For fixed $(i,j)\in\{ 1,\ldots,I\}\times\{ 1,\ldots,J\}$ and strategies 
$(\alpha,\beta)\in\AR(t)\times\BR(t)$, the payoff of the
game with only one possible terminal payoff function $g_{ij}$ will be denoted by
\[ J_{ij}(t,x,\alpha,\beta)=E[g_{ij}(X^{t,x,\alpha,\beta}_T)].\]
Now let $(\overline\alpha,\overline\beta)\in \AR_r(t)\times\BR_r(t)$ be two random strategies, with
$\overline\alpha=(\alpha^1,\ldots,\alpha^R;r^1,\ldots,r^R)$ and
$\overline\beta=(\beta^1,\ldots,\beta^S;s^1,\ldots,s^S)$.
The payoff  associated with the pair  $(\overline\alpha,\overline\beta)\in \AR_r(t)\times\BR_r(t)$),
is the average of the payoffs with respect to the probability distributions 
associated to the strategies:
\[J_{ij}(t,x,\overline\alpha,\overline\beta)=
\sum_{k=1}^R\sum_{l=1}^Sr^ks^lE[g_{ij}(X^{t,x,\alpha^k,\beta^l}_T)].\]
Further, for $p\in\Delta(I)$, $j\in\{ 1,\ldots,J\}$, 
$\hat\alpha\in(\AR_r(t))^I$ and $\overline\beta\in\BR_r(t)$ 
 we will use the notation
\[ J^p_j(t,x,\hat\alpha,\overline\beta)=
\sum_{i=1}^Ip_iJ_{ij}(t,x,\overline\alpha_i,\overline\beta)
=\sum_{i=1}^Ip_i\sum_{k,l}r^ks^lE[g_{ij}(X^{t,x,\alpha^k_i,\beta^l}_T)].\]\;
A symmetric notation holds for $\overline\alpha\in\AR_r(t)$ and
 $\hat\beta\in(\BR_r(t))^J$. Finally, the payoff of the game is, for 
$(\hat\alpha,\hat\beta)\in(\AR_r(t))^I\times(\BR_r(t))^J$, $p\in \Delta(I)$, $q\in \Delta(J)$,
\[ J^{p,q}(t,x,\hat\alpha,\hat\beta)=\sum_{i=1}^I\sum_{j=1}^Jp_iq_j
J_{ij}(t,x,\overline\alpha_i,\overline\beta_j).\]
 The reference to $(t,x)$ in the notations is dropped when there is no possible confusion : we 
will write $J_{ij}(\alpha,\beta), J_{ij}(\overline\alpha,\beta),\ldots$.
\\

We define the value functions for the game by
\[\begin{array}{c}
V^+(t,x,p,q)=\inf_{\hat\alpha\in(\AR_r(t))^I}\sup_{\hat\beta\in
\BR_r(t))^J}J^{p,q}(t,x,\hat\alpha,\hat\beta),\\
V^-(t,x,p,q)=\sup_{\hat\beta\in
(\BR_r(t))^J}\inf_{\hat\alpha\in(\AR_r(t))^I}J^{p,q}(t,x,\hat\alpha,\hat\beta).
\end{array}\]
Again we will write  $V^+(p,q)$ and $V^-(p,q)$ if there is no possible confusion on $(t,x)$.\\

 The following lemma follows easily from classical estimations for stochastic 
differential equations :

\begin{Lemma}\label{reguV+V-}
$V^+$ and $V^-$ are bounded, Lipschitz continuous with respect to $x,p,q$ and H\"older continuous with
respect to $t$.
\end{Lemma}

Following \cite{AM} we now state one of the basic properties of the value functions. The technique of proof of this statement is
known as the splitting method in repeated game theory (see \cite{AM}, \cite{So}).

\begin{Proposition}\label{CvV+V-}
For all $(t,x)\in[0,T]\times\R^n$, the maps $(p,q)\rightarrow V^+(t,x,p,q)$ and
$(p,q)\rightarrow V^-(t,x,p,q)$ are convex in $p$ and concave in $q$.
\end{Proposition}

\noindent{\em Proof:}
We only prove the result for $V^+$,  the proof for $V^-$ is the same.
First $V^+$ can be rewritten as
\[ V^+(p,q)=\inf_{\hat\alpha\in(\AR_r(t))^I}\sum_{j=1}^Jq_j
\sup_{\overline\beta\in\BR_r(t)}J^p_j(\hat\alpha,\overline\beta).\]
It follows that $V^+$ is concave in $q$.\\

Now fix $q\in\Delta(J)$ and let $p,p'\in\Delta(I)$
and $a\in(0,1)$. Without loss of generality we can assume that, for all $i\in\{1,\ldots,I\}$, $p_i$ and $p'_i$ 
are not simultaneously equal to zero.\\
We get a new element of $\Delta(I)$ if we set
 $p^a=ap+(1-a)p'$.
For $\epsilon>0$, let $\hat\alpha\in(\AR_r(t))^I$ be $\epsilon$-optimal
for $V^+(p,q)$ (resp. $\hat\alpha'\in(\AR_r(t))^I$ $\epsilon$-optimal for $V^+(p',q)$).\\
We define a new strategy 
$\hat\alpha^a=(\overline\alpha_1^a,\ldots,\overline\alpha_I^a)$ by
\[ \overline\alpha_i^a=(\alpha_i^1,\ldots,\alpha_i^R,
\alpha'^1_i,\ldots,\alpha'^{R'}_i;
(r_i^a)^{1},\ldots, (r_i^a)^{(R+R')}), i\in\{ 1,\ldots, I\},\]
with
\[(r_i^a)^{k}=\left\{
\begin{array}{ll}
\frac{ap_i}{p^a_i}r_{i}^{k} &\mbox{for } k\in\{ 1,\ldots ,R\},\\
\\
\frac{(1-a)p_i'}{p^a_i}r'^{k-R}_i &\mbox{for } k\in\{ R+1,\ldots ,R+R'\}
\end{array}\right.\]
(it is easy to check that $\hat\alpha^a\in(\AR_r(t))^I$).\\
This means that, for all $\hat\beta\in(\BR_r(t))^J$,
\[
J^{p^a,q}(\hat\alpha^a,\hat\beta)=
\sum_{i=1}^I\big\{
ap_i \sum_{k=1}^{R}r_i^k
J_i^q(\alpha_i^k,\hat\beta)
+(1-a)p_i' \sum_{k=1}^{R'}r'^k_i
J_i^q(\alpha'^k_i,\hat\beta)\big\}
\]

Thus
\[
\sup_{\hat\beta\in(\br(t))^J}J^{p^a,q}(\hat\alpha^a,\hat\beta)
\leq a\sup_{\hat\beta\in(\BR_r(t))^J}J^{p,q}(\hat\alpha,\hat\beta)
+(1-a)\sup_{\hat\beta\in(\BR_r(t))^J}J^{p',q}(\hat\alpha',\hat\beta).
\]
It follows by the choice of $\hat\alpha$ and $\hat\alpha'$ that
\[ V^+(p^a,q)\leq aV^+(p,q)+(1-a)V^+(p',q).\]
\cq

\section{Subdynamic programming and Hamilton-Jacobi-\allowbreak Bellman equations 
for the Fenchel conjugates.}
\label{subdynamic}

Since $V^+$ and $V^-$ are convex with respect to $p$ and concave with respect to $q$, it is natural to introduce the
Fenchel conjugates of these functions. For this we use the following notations.\\
For any $w:[0,T]\times\R^n\times\Delta(I)\times\Delta(J)\rightarrow\R$,
we define the Fenchel conjugate $w^*$ of $w$ with respect to $p$ by
\[ w^*(t,x,\hat p,q)=\sup_{p\in\Delta(I)}\{ \langle \hat p, p\rangle -w(t,x,p,q)\},
\; (t,x,\hat p,q)\in [0,T]\times\R^n\times\R^I\times\Delta(J).\]
For $w$ defined on the dual space $[0,T]\times\R^n\times\R^I\times\Delta(J)$, we also set
\[ w^*(t,x,p,q)=\sup_{\hat p\in\R^I}\{ \langle \hat p, p\rangle -w(t,x,\hat p,q)\},
\; (t,x,p,q)\in [0,T]\times\R^n\times\Delta(I)\times\Delta(J).\]
It is well known that, if $w$ is convex, we have $(w^*)^*=w$.\\
We also have to introduce the concave conjugate with respect to $q$ of a map $w:[0,T]\times\R^n\times\Delta(I)\times\Delta(J)\rightarrow\R$:
\[ w^\sharp(t,x,p,\hat q)=\inf_{q\in\Delta(J)}\{ \langle \hat q, q\rangle -w(t,x,p,q)\},
\; (t,x,p,\hat q)\in[0,T]\times\R^n\times\Delta(I)\times\R^J.\]
We use the following notations for the sub- and superdifferentials with respect to $\hat p$ and $\hat q$
respectively: if $w: [0,T]\times \R^n\times \R^I\times \Delta(J)\to \R$,  we set
$$
\partial^-_{\hat p} w(t,x,\hat p, q)=\{ p \in\R^I, \; w(t,x,\hat p, q)+\langle p, \hat p'-\hat p\rangle \leq w(t,x,\hat p', q)\; \forall \hat p'\in \R^I\}
$$
and if $w: [0,T]\times \R^n\times \Delta(I)\times \R^J\to \R$
$$
\partial^+_{\hat q} w(t,x,p, \hat q)=\{ q \in\R^J, \; w(t,x,p, \hat q)+\langle q, \hat q'-\hat q\rangle  \geq w(t,x,p, \hat q')\; \forall \hat q'\in \R^J\}.
$$

In this chapter, we will show that $V^{+\sharp}$ and $\vme$ satisfy
a subdynamic programming property. This part follows several ideas of \cite{DM2}, \cite{DM1}.

\begin{Lemma}(Reformulation of $\vme$)\\
For all $(t,x,\hat p,q)\in[0,T]\times\R^n\times\R^I\times\Delta(J)$, we have
\begin{equation}
\label{refv}
\vme(t,x,\hat p,q)=
\inf_{\hat\beta\in(\BR_r(t))^J}\sup_{\alpha\in\AR(t)}\max_{i\in\{ 1,\ldots,I\}}
\Big\{ \hat p_i-J^q_i(t,x,\alpha,\hat\beta)\Big\}.
\end{equation}
\end{Lemma}

\noindent{\em Proof.} 
We begin to establish a first expression for $\vme$:
\begin{equation}
\label{refvr}
\vme(\hat p,q)=
\inf_{\hat\beta\in(\BR_r(t))^J}
\sup_{\overline\alpha\in\ar(t)}\max_{i\in\{ 1,\ldots,I\}}
\Big\{ \hat p_i-J^q_i(\overline\alpha,\hat\beta)\Big\}
\end{equation}
(the difference with (\ref{refv})
is that player I here can use random strategies.)\\
Let's denote by $e=e(\hat p,q)$ the right hand term of (\ref{refvr}). 
First we prove that $e$ is convex with respect to $\hat p$ :\\
Fix $q\in\Delta(J)$, $\hat p,\hat p'\in\R^I$ and  $a\in(0,1)$.\\
For $\epsilon>0$, let $\hat\beta$  (resp. $\hat\beta'$)$\in(\br(t))^J$
be some
$\epsilon$-optimal strategy for $e(\hat p,q)$ (resp. $e(\hat p',q)$).\\
Set $\hat p^a=a\hat p+(1-a) \hat p'$.\\
We define a new strategy $\hat\beta^a \in(\BR_r(t))^J$ by
\[ \overline\beta^a_j=(\beta_j^1,\ldots, \beta_j^S,
\beta'^1_j,\ldots, \beta'^{S'}_j;(s_j^a)^1,\ldots,(s_j^a)^{S+S'}),\;\;
j\in\{ 1,\ldots,J\},\]
with
\[ (s^a)^k_j=\left\{\begin{array}{ll}
as_j^k &\mbox{for } k\in\{ 1,\ldots, S\},\\
\\
(1-a)s'^{k-S}_j& k\in\{ S+1,\ldots, S+S'\}.
\end{array}\right.\]
Let $\overline\alpha\in\AR_r(t)$. 
Since the application $(x_1,\ldots,x_I)\rightarrow\max\{x_i,i=1,\ldots,I\}$ is convex, 
we have
\[ \begin{array}{rl}
\max_i\Big\{ \hat p^a_i-J^q_i(\overline\alpha,\hat\beta^a)\Big\}=&
\max_i\Big\{ a(\hat p_i-J^q_i(\overline\alpha,\hat\beta))
+(1-a)(\hat p'_i-J^{q'}_i(\overline\alpha,\hat\beta'))\Big\}\\
\\
\leq &a\sup_{\overline\alpha\in\AR_r(t)}
\max_i (\hat p_i-J^q_i(\overline\alpha,\hat\beta^a))\\
&+(1-a)\sup_{\overline\alpha\in\AR_r(t)}
\max_i(\hat p^1_i-J^{q'}_i(\overline\alpha,\hat\beta'))\\
\\
\leq & ae(\hat p,q)+(1-a)e(\hat p',q)+\epsilon.
\end{array}\]
Since $\epsilon$ is arbitrary, we can deduce that $e$ is convex with respect to $\hat p$.\\
The next step is to prove that $e^*=V^-$. By the convexity of $e$, this will imply that
$V^{-*}=e$.\\
We can reorganize $e^*(p,q)$ as follows :
\[
\begin{array}{rl}
e^*(p,q)=&\sup_{\hat p\in\R^I}\left\{ \sum_{i=1}^I\hat p_ip_i+\sup_{\hat\beta\in(\BR_r(t))^J}
\inf_{\overline\alpha\in\AR_r(t)}\min_{i'\in\{ 1,\ldots, I\}}
\{ J_{i'}^q(\overline\alpha,\hat\beta)-\hat p_{i'}\}\right\}\\
\\
=&\sup_{\hat\beta\in(\BR_r(t))^J}\sup_{\hat p\in\R^I}\sum_{i=1}^I
p_i\min_{i'\in\{ 1,\ldots, I\}}\left\{\inf_{\overline\alpha\in\AR_r(t)}
 J_{i'}^q(\overline\alpha,\hat\beta)+(\hat p_i-\hat p_{i'})\right\}
\end{array}\]
The supremum over $\hat p\in\R^I$ is attained for $\hat p_{i'}=\inf_{\overline\alpha\in\AR_r(t)}
 J_{i'}^q(\overline\alpha,\hat\beta)$ and we get the claimed result.\\

Finally, to get (\ref{refv}),
 it remains to show that player I can use non random strategies.\\
Indeed, writing $\vme$ as in (\ref{refvr}) and since $\AR(t)\subset\ar(t)$, 
it is obvious that the left hand side of (\ref{refv})
is not smaller than the right hand side.\\
Concerning the reverse inequality, we can write
\[\begin{array}{rl}
\sup_{\overline\alpha\in\ar(t)}&
\max_i\Big\{ \hat p_i-J^q_i(\overline\alpha,\overline\beta)\Big\}\\
\\
&\leq\sup_{R\in\N^*}\sup_{(\alpha^1,\ldots,\alpha^R)\in(\AR(t))^R,
(r^1,\ldots,r^R)\in\Delta(R)}
\sum_{k=1}^Rr^k\max_i\Big\{\hat p_i-J^q_i(\alpha^k,\hat\beta)\Big\}\\
\\
&\leq\sup_{R\in\N^*}\sup_{
(r^1,\ldots,r^R)\in\Delta(R)}
\sum_kr^k\sup_{\alpha\in\AR(t)}\max_i\Big\{
\hat p_i-J^q_i(\alpha,\hat\beta)\Big\}.
\end{array}
\]
The result follows after one recalls that $\sum_{k=1}^Rr^k=1$.\\

\cq

\begin{Proposition}
\label{dynprog}
(Subdynamic programming for $V^{-*}$)\\
For all $0\leq t_0\leq t_1\leq T, x_0\in\R^n, \hat p\in\R^I,q\in\Delta(J)$, 
it holds that
\[ V^{-*}(t_0,x_0,\hat p,q)\leq\inf_{\beta\in\BR(t_0)}\sup_{\alpha\in\AR(t_0)}
E[V^{-*}(t_1,X_{t_1}^{t_0,x_0,\alpha,\beta},\hat p,q)].\]
\end{Proposition}

\dem\ Set $\vome(t_0,t_1,x_0,\hat p,q)=
\inf_{\beta\in\BR(t_0)}\sup_{\alpha\in\AR(t_0)}
E[\vme(t_1,X_{t_1}^{t_0,x_0,\alpha,\beta},\hat p,
q)]$. \\
For $\epsilon>0$, let $\beta^{\epsilon}\in\BR(t_0)$ be $\epsilon$-optimal for 
$\vome(t_0,t_1,x_0,\hat p,q)$, and, for all $x\in\R^n$, let $\hat\beta^x\in(\br(t_1))^J$
be $\epsilon$-optimal for $\vme(t_1,x,\hat p,q)$. 
By the uniformly Lipschitz assumptions 
for the parameters of the dynamics, there exists $R>0$
such that, for all $\alpha\in\AR(t_0)$,
\[ P[X_{t_1}^{t_0,x_0,\alpha,\beta^\epsilon}\in B(x_0,R)]
\geq 1-\epsilon,\]
where $B(x_0,R)$ denotes the ball in $\R^n$ of center $x_0$ and radius $R$.\\
Remark that $J^q_i$ and $V^{-*}$ are uniformly Lipschitz continuous in $x$.
This implies that we can find $r>0$ such that, for any $x\in\R^n$ and $y\in\BR(x,r)$, 
$\hat\beta^x$ is $2\epsilon$-optimal for $V^{-*}(t_1,y,\hat p,q)$.\\
Now let $x_1,\ldots,x_M\in\R^n$ such that $\cup_{m=1}^MB(x_m,\frac r2)
\supset B(x_0,R)$.\\
Set $\hat\beta^m=\hat\beta^{x_m}$ for $m=1,\ldots,M$ and choose some arbitrary
$\hat\beta^0\in(\br(t_1))^J$.\\
Each $\hat\beta^m$ is detailed in the following way:
\[ \hat\beta^m=(\overline\beta^m_1,\ldots,\overline\beta_J^m),\]
with
\[\overline\beta^m_j=
(\beta_j^{m,1},\ldots,\beta_j^{m,S^m_j};s^{m,1}_j,\ldots,s_j^{m,S^m_j}).\]

Let $\delta$ be a common delay for $\hat\beta^0,\ldots,\hat\beta^M$
 that we can choose as small as we need :\\
$0<\delta<\frac{r^2\epsilon}{4C}\wedge (t_1-t_0)$, where $C>0$ is defined through
the parameters of the dynamics by
\[ \forall\alpha\in\AR(t),\beta\in\BR(t),t,t'\in[t_0,T],\;
E[|X_t^{t_0,x_0,\alpha,\beta}-X_{t'}^{t_0,x_0,\alpha,\beta}|^2]\leq C|t-t'|. 
\]
We then have in particular, for all $\alpha\in\AR(t)$ and $\beta\in\BR(t)$,
\begin{equation}
\label{delta}
P[|X_{t_1}^{t_0,x_0,\alpha,\beta}-X_{t_1-\delta}^{t_0,x_0,\alpha,\beta}|>\frac r2]\leq\epsilon.
\end{equation}
Let $(E_m)_{m=1,\ldots,M}$ be a Borel measurable partition of $B(x_0,R)$, 
such that, for all $m\in\{ 1,\ldots, M\}$, $E_m\subset B(x_m,\frac r2)$. 
Set $E_0=B(x_0,R)^c$.\\
We are now able to 
define a new strategy for player II, $\hat\beta^\epsilon\in
(\br(t_0))^J$:\\
Fix $j\in\{ 1, \ldots,J\}$. For $l=(l_0,\ldots,l_M)\in L:=
\Pi_{m=0}^M\{ 1,\ldots,S^m_j\}$,
set $s_j^l=\Pi_{m=0}^Ms^{m,l_m}_j$.
Remark that $\{ s^l_j,l\in L\}\in\Delta(\mbox{Card}(L))$.\\
Then, for $l\in L, l=(l_0,\ldots ,l_M)$, we define $(\beta_j^{\epsilon})^l
\in\BR(t_0)$ by
\[ \begin{array}{rl}
\forall f\in C([t_0,T],\R^n),\forall t\in[t_0,T],\\
(\beta_j^\epsilon)^l(t,f)&=\left\{
\begin{array}{ll}
\beta^\epsilon(t,f)& \mbox{ if } t\in [t_0,t_1),\\
\beta^{m,l_m}_j(t,f|_{[t_1,T]})&\mbox{ if } t\in [t_1,T] \mbox{ and }
 f(t_1-\delta)\in E_m.
\end{array}\right.\end{array}\]
We set $\overline\beta_j^\epsilon:=((\beta_j^\epsilon)^l;s_j^l,l\in L)\in\br(t_0)$, 
and finally
$\hat\beta^\epsilon=(\overline\beta_1^\epsilon,\ldots, \overline\beta_J^\epsilon)$.\\

\noindent For some fixed $\alpha\in\AR(t_0)$ and $f\in C([t_0,t_1],\R^n)$, we
define a new strategy $\alpha_f\in\AR(t_1)$ by:\\
 for all $t\in[0,T]$ and $f'\in C([t_1,T],\R^n)$,
\[ \alpha_f(t,f')=\alpha(t,\tilde f),
\mbox{ with }
\tilde f(t)=\left\{
\begin{array}{l}f(t)\mbox{ for } t\in [t_0,t_1],\\
f'(t)-f'(t_1)+f(t_1),\mbox{ for } t\in (t_1,T].
\end{array}\right.\]

Set
$X_\cdot^\epsilon=
X_\cdot^{t_0,x_0,\alpha,\beta^{\epsilon}}$ and, 
for $m\in\{ 0,\ldots,M\}$,  $A_m=\{ X_{t_1-\delta}^\epsilon\in E_m\}$.
Set further
$A=\{ |X_{t_1}^\epsilon-X_{t_1-\delta}^\epsilon|\leq\frac r2\}$.
By (\ref{delta}), it holds that $P[A^c]\leq\epsilon$.
Remark also that, on each $A\cap A_m$, $X^\epsilon_{t_1}$ belongs to $B(x_m,r)$ and consequently,
still on $A\cap A_m$, $\hat\beta^m$ is $2\epsilon$-optimal for $V^{-*}(t_1,X_{t_1}^\epsilon,\hat p,q)$.\\
For all $i\in\{ 1,\ldots, I\},j\in\{ 1, \ldots ,J\}$
and $l\in L$,
we have
\[
E[g_{ij}(X_T^{t_0,x_0,\alpha,(\beta_j^\epsilon)^l})|\FR_{t_1}]=
\sum_{m=0}^M\i_{A_m}
E[g_{ij}(X_T^{t_1,y,
\alpha_f,\beta^{m,l_m}_j})]|_{y=X_{t_1}^\epsilon, f=X_\cdot^\epsilon|_{[t_0,t_1]}}.
\]
It follows that
\[\begin{array}{rl}
J^q_i(t_0,x_0,\alpha,\hat\beta^\epsilon)=&
\sum_{j=1}^Jq_j\sum_{l\in L}s_j^lE[g_{ij}(X_T^{t_0,x_0,\alpha,(\beta^\epsilon_j)^l})]\\
\\
=& E[\sum_{m=0}^M\i_{A_m}J^q_i(t_1,X^\epsilon_{t_1},\alpha_{X^\epsilon|_{[t_0,t_1]}},\hat\beta^m)].
\end{array}\]
And
\[\begin{array}{rl}
\max_{i\in\{ 1,\ldots,I\}}\Big\{ \hat p_i-&
J^q_i(t_0,x_0,\alpha,\hat\beta^\epsilon)
\Big\}\\
\leq&
E[\sum_{m=0}^M\i_{A_m}\max_{i\in\{ 1,\ldots,I\}}
\{ \hat p_i-J_i^q(t_1,X_{t_1}^\epsilon,\alpha_{X_\cdot^\epsilon|_{[t_0,t_1]}},
\hat\beta^m)\}]\\
\\
\leq &
E[\sum_{m=0}^M\i_{A_m}(\sup_{\alpha\in\AR(t_1)}
\max_{i\in\{ 1,\ldots,I\}}\{ \hat p_i-J_i^q(t_1,X_{t_1}^\epsilon,\alpha,\hat\beta^m)\})]\\
\\
\leq &
E[(V^{-*}(t_1,X_{t_1}^\epsilon,\hat p,q)+2\epsilon)\i_{A\cap\{ X^\epsilon_{t_1}\in B(x_0,R)\}}]\\
&\hspace{3cm}+\max_{i\in\{ 1,\ldots,I\}}\{ |\hat p_i|+K\}(P[A^c]+P[X^\epsilon_{t_1}\not\in B(x_0,R)]),
\end{array}\]
by the choice of $(\hat\beta^m,m\in\{ 1,\ldots,M\})$ and where $K$ is an upper bound of $|g|$.\\
By the choice of $R$ and with the notation $K(\hat p)=4\max_{i\in\{ 1,\ldots,I\}}\{ |\hat p_i|+K\}
+\epsilon$,
we get
\[\begin{array}{rl}
\max_{i\in\{ 1,\ldots,I\}}\Big\{ \hat p_i-
J^q_i(t_0,x_0,\alpha,\hat\beta^\epsilon)
\Big\}
\leq&
E[V^{-*}(t_1,X_{t_1}^\epsilon,\hat p,q)+2\epsilon]+K(\hat p)\epsilon\\
\\
\leq &
\sup_{\alpha\in\AR(t_0)}E[V^{-*}(t_1,X_{t_1}^{t_0,x_0,\alpha,\beta^\epsilon},\hat p,q)]+2\epsilon(1+K(\hat p))\\
\\
\leq &
V^{-*}_1(t_0,t_1,x_0,\hat p,q)+\epsilon(3+2K(\hat p))
\end{array}\]
(for the last inequality, recall that
$\beta^\epsilon$ was chosen $\epsilon$-optimal for $\vme_1(t_0,t_1,x_0,\hat p,q)$).\\
We can deduce the result.
\cq

A classical consequence of the subdynamic programming principle for $V^{-*}$ is that this function
is a subsolution of some associated Hamilton-Jacobi equation. We give a proof of that result for sake
of completeness. 

\begin{Corollary}\label{HJIV-} For any $(\hat p, q)\in \R^I\times \Delta(J)$, 
$\vme(\cdot,\cdot, \hat p, q)$ is a subsolution in the viscosity sense of
\[ w_t+H^{-*}(t,x,Dw,D^2w)=0,\qquad  (t,x)\in(0,T)\times \R^n,\]
with 
\begin{equation}\label{H--}\begin{array}{l}
H^{-*}(t,x,p,A)=-H^-(t,x,-p,-A)=\\
\qquad \inf_{v\in V}\sup_{u\in U}\{ \langle b(t,x,u,v),p\rangle +\frac 12\mbox{Tr}(A\sigma(t,x,u,v)
\sigma^*(t,x,u,v))\}.
\end{array}
\end{equation}
\end{Corollary}

\dem\
For $(t_0,x_0)\in[0,T]\times\R^n, \hat p\in\R^I, q\in\Delta(J)$ fixed, let
$\phi\in C^{1,2}$ such that
$\phi(t_0,x_0)=\vme(t_0,x_0,\hat p,q)$ and, for all $(s,y)\in[0,T]\times\R^n$, 
$\phi(s,y)\geq \vme(s,y,\hat p,q)$.
\\
We have to prove that
\[ \phi_t(t_0,x_0)+H^{-*}(t_0,x_0,D\phi(t_0,x_0),D^2(t_0,x_0))\geq 0.\]
Suppose that this is false and consider $\theta>0$ such that
\begin{equation}
\label{theta}
\phi_t(t_0,x_0)+H^{-*}(t_0,x_0,D\phi(t_0,x_0),D^2(t_0,x_0))\leq -\theta<0.
\end{equation}
Set $\Lambda(t,x,u,v)=\phi_t(t,x)+\langle b(t,x,u,v),D\phi(t,x)\rangle
+\mbox{Tr}(D^2\phi(t,x)\sigma(t,x,u,v)
\sigma^*(t,x,u,v))$.
Since, for fixed $\hat p$, $V^{-*}$ is bounded, we can choose $\phi$ such that $\phi_t$ and $D^2\phi$ are 
also bounded. It follows that, for some $K>0$,
we have $|\Lambda(t,x,u,v)|\leq K$.
\\
Now the relation (\ref{theta}) is equivalent to
\[
 \inf_{v\in V}\sup_{u\in U}\Lambda(t_0,x_0,u,v)\leq -\theta\;.
\]
This implies the existence of a control  $v_0\in V$ such that, for all $u\in U$,
\[
\Lambda(t_0,x_0,u,v_0)\leq-\frac{2\theta}3.
\]
Moreover, since $\Lambda$ is continuous in $(t,x)$, uniformly in $u,v$,
we can find $R>0$ such that, 
\begin{equation}
\label{8}
\forall(t,x)\in[t_0,T]\times\R^n,
|t-t_0|\vee\|x-x_0\|<R, \forall u\in U,
 \Lambda(t,x,u,v_0)\leq-\frac\theta 2.
\end{equation}
Now define a strategy for player II by $\beta_0(t,f)=v_0$ for all 
$(t,f)\in[t_0,T]\times C([t_0,T],\R^n)$.\\
Fix $\epsilon>0$ and $t\in (t_0,R)$.
Because of the subdynamical programming (Proposition \ref{dynprog}), 
there exists $\alpha_{\epsilon,t}\in\AR(t_0)$
such that
\begin{equation}
\label{utidyn}
E[V^{-*}(t_1,X_{t_1}^{t_0,x_0,\alpha_{\epsilon,t},\beta_0},\hat p,q)]-V^{-*}(t_0,x_0,\hat p,q)\geq -\epsilon(t-t_0).
\end{equation}
Let $(u_s,v_s)\in\UR(t_0)\times\VR(t_0)$ the controls associated to $(\alpha_{\epsilon,t},\beta_0)$ 
by the relation (\ref{alphau})
and set $X_\cdot=X_\cdot^{t_0,x_0,\alpha_{\epsilon,t},\beta_0}=X_\cdot^{t_0,x_0,u,v}$. 
(Remark that, by the choice of $\beta_0$, 
$(v_s)$ is constant and equal to $v_0$.)\\
Now we write It\^o's formula for $\phi(t,X_t)$:
\begin{equation}
\label{itof}
\begin{array}{rl}
\phi(t,X_t)-\phi(t_0,x_0)=&
\int_{t_0}^t\Lambda(s,X_s,u_s,v_s)ds\\
&+\int_{t_0}^t\langle D\phi(s,X_s),b(s,X_s,u_s,v_s)\rangle dB_s.
\end{array}
\end{equation}
By (\ref{utidyn}), (\ref{itof}) and the definition of $\phi$, we have
\begin{equation}
\label{epsit}
E[\int_{t_0}^t\Lambda(s,X_s,u_s,v_s)ds]\geq -\epsilon(t-t_0).
\end{equation}
In the other hand, there exists a constant $C>0$ depending only on the parameters of $X$, such that
\[
 P[\| X_\cdot-x_0\|_t> R]\leq \frac{C(t-t_0)^2}{R^4},
\]
with the notation $\|f\|_t=\sup_{s\in[t_0,t]}\| f(s)\|$.\\
Following (\ref{8}), this implies that, for all $t\in [t_0,T\wedge (t_0+R)]$,
\begin{equation}
\label{majo}
E\left[ I_{\{ \| X_\cdot-x_0\|_t<R\}}\int_{t_0}^t\Lambda (s,X_s, u_s,v_s)ds
 \right]\leq-\frac\theta 2(t-t_0).
\end{equation}
By (\ref{epsit}) and (\ref{majo}), we now have  
\[\begin{array}{rl}
-\epsilon(t-t_0)\leq &
E[\int_{t_0}^t\Lambda(s,X_s,u_s,v_s)dsI_{\{ \| X_\cdot-x_0\|_t>R\}}]
+
E[\int_{t_0}^t\Lambda(s,X_s,u_s,v_s)dsI_{\{ \| X_\cdot-x_0\|_t\leq R\}}]\\
\\
\leq &
\frac{KC}{R^4}(t-t_0)^2-\frac\theta 2(t-t_0),
\end{array}\]
or, equivalently,
\[ \frac\theta 2\leq\frac{KC}{R^4}(t-t_0)+\epsilon.
\]
Since $t-t_0$ and $\epsilon$ can be chosen arbitrarily small, we get a contradiction.\cq

For $V^+$ we have:

\begin{prop}(Superdynamic programming and HJI equation for $V^{+\sharp}$)\label{CasV+} \\
For all $0\leq t_0\leq t_1\leq T, x_0\in\R^n, p\in\Delta(I),\hat q\in\R^J$, 
it holds that
\[ V^{+\sharp}(t_0,x_0,p,\hat q)\geq\inf_{\beta\in\BR(t_0)}\sup_{\alpha\in\AR(t_0)}
E[V^{+\sharp}(t_1,X_{t_1}^{t_0,x_0,\alpha,\beta},p,\hat q)].\]
As a consequence, for any $(p,\hat q)\in \Delta(I)\times \R^J$, $V^{+\sharp}(\cdot,\cdot,p,\hat q)$ is a supersolution in viscosity sense of
\[ w_t+H^{+*}(t,x,Dw,D^2w))=0, \qquad (t,x)\in(0,T)\times \R^n,\]
where 
\begin{equation}\label{H++}\begin{array}{l}
H^{+*}(t,x,p,A)=-H^+(t,x,-p,-A)=\\
\qquad \sup_{u\in U}\inf_{v\in V}\{ \langle b(t,x,u,v),p\rangle +\frac 12\mbox{Tr}(A\sigma(t,x,u,v)
\sigma^*(t,x,u,v))\}.
\end{array}
\end{equation}

\end{prop}

{\bf Proof : } We note that $V^+$ is equal to the opposite of the lower  value of the game in which we replace $g_{ij}$ by $-g_{ij}$,
Player I is the maximizer and in which the respective roles of $p$ and $q$ are exchanged. Using Proposition \ref{dynprog}
in this framework gives the superdynamic programming principle. Now Corollary \ref{HJIV-} shows that, for any 
$(p,\hat q)\in \Delta(I)\times \R^J$, $(-V^+)^*(\cdot, \cdot, p, \hat q)=-V^{+\sharp}(\cdot, \cdot, p, -\hat q)$ is a subsolution of 
\[ w_t+H^+(t,x,Dw,D^2w))=0, \qquad (t,x)\in(0,T)\times \R^n.\]
Hence $V^{+\sharp}(\cdot, \cdot, p, -\hat q)$ is a supersolution of 
\[ w_t+H^{+*}(t,x,Dw,D^2w))=0,  \qquad (t,x)\in(0,T)\times \R^n.\]
Since this holds true for any $(p,\hat q)$, this proves our claim.
\cq


\section{Comparison principle and existence of a value}

In this section we first state a new comparison principle and apply it to get the existence and the characterization of the value.
Then we give a proof for the comparison principle.

\subsection{Statement of the comparison principle and existence of a value}

Let $H: [0,T]\times \R^n\times \R^n\times {\cal S}_n\times \Delta(I)\times \Delta(J)\to \R$ be continuous and satisfy
\begin{equation}\label{StructH}
\begin{array}{l}
H(s,y,\xi_2,X_2,p,q)-H(t,x,\xi_1,X_1,p,q) \; \geq \\
\qquad  - \omega\left(|\xi_1-\xi_2|+a |(t,x)-(s,y)|^2+b+ |(t,x)-(s,y)|(1+|\xi_1|+|\xi_2|)\right)
\end{array}
\end{equation}
where $\omega$ is continuous and non decreasing with $\omega(0)=0$, 
for any $a,b\geq 0$, $(p,q)\in \Delta(I)\times\Delta(J)$, $s,t\in [0,T]$, $x,y,\xi_1,x_2\in \R^n$ and $X_1,X_2\in {\cal S}_n$ such that 
$$
 \left(\begin{array}{ll} -X_1 & 0\\ 0 & X_2\end{array}\right) \leq a\left(\begin{array}{cc}I&-I\\-I& I\end{array}\right)+bI
$$

\begin{Definition}
We say that a  map $w:(0,T)\times \R^n\times \Delta(I)\times \Delta(J)\to \R$
is a {\rm supersolution in the dual sense} of equation
\begin{equation}\label{HJI}
w_t+H(t,x,Dw,D^2w,p,q)=0
\end{equation}
if  $w=w(t,x,p,q)$ is lower semicontinuous, concave with respect to $q$ and if, for any ${\cal C}^2((0,T)\times \R^n)$ function $\phi$ such that $(t,x)\to w^*(t,x,\hp,\barq)-\phi(t,x)$ has
a maximum at some point  $(\bart,\barx)$ for some $(\hp,\barq)\in \R^I\times\Delta(J)$, we have
$$
\phi_t(\bart,\barx)-H(\bart, \barx, -D\phi(\bart, \barx), -D^2\phi(\bart, \barx), p, \barq)\geq 0
\qquad \forall p \in \partial^-_\hp w^*(\bart,\barx, \hp,\barq)\;.
$$
We say that $w$ is a {\rm subsolution of (\ref{HJI}) in the dual sense} if
$w$ is upper semicontinuous, convex with respect to $p$ and if, for any ${\cal C}^2((0,T)\times \R^n)$ function $\phi$ such that $(t,x)\to w^\sharp(t,x,\barp,\hq)-\phi(t,x)$ has
a minimum at some point  $(\bart,\barx)$ for some $(\barp,\hq)\in \Delta(I)\times \R^J$, we have
$$
\phi_t(\bart,\barx)-H(\bart, \barx, -D\phi(\bart, \barx), -D^2\phi(\bart, \barx), \barp, q)\leq 0
\qquad \forall q \in \partial^+_\hq w^\sharp(\bart,\barx, \barp, \hq)\;.
$$
A solution of (\ref{HJI}) in the dual sense is a map which is sub- and supersolution in the dual sense.
\end{Definition}

{\bf Remarks : } \begin{enumerate}
\item We have proved in Corollary \ref{HJIV-} that $V^-$ is a dual supersolution of the HJ equation
$$
w_t+H^-(t,x,Dw,D^2w)=0\;,
$$
where $H^-$ is defined by (\ref{H--}), while Proposition \ref{CasV+}  shows that $V^+$ is a dual subsolution of the HJ equation
$$
w_t+H^+(t,x,Dw,D^2w)=0\;,
$$
where $H^-$ is defined by (\ref{H++}).

\item The necessity to deal with a Hamiltonian $H$ with a $(p,q)$ dependence will become clear in the next section where we study
differential games with running costs.

\item An equivalent definition of the notion of dual super- or subsolution in given in Lemma \ref{CaracSuperSol} below.
\end{enumerate}

The main result of this section is the following:

\begin{thm}[Comparison principle]\label{comparison} Let us assume that $H$
satisfies the structure condition (\ref{StructH}).
Let $w_1$ be a bounded, H\"{o}lder continuous subsolution of (\ref{HJI}) in the dual sense which is uniformly Lipschitz continuous w.r. to $q$
and  $w_2$ be a bounded, H\"{o}lder continuous  supersolution of (\ref{HJI}) in the dual sense which is uniformly Lipschitz continuous w.r. to $p$. Assume that 
\begin{equation}\label{InitCond}
w_1(T,x,p,q)\leq  w_2(T,x,p,q)\qquad \forall (x,p,q)\in \R^n\times \Delta(I)\times \Delta(J)\;.
\end{equation}
Then
$$
w_1(t,x,p,q)\leq w_2(t,x,p,q)\qquad \forall (t,x,p,q)\in [0,T]\times \R^n\times \Delta(I)\times \Delta(J)\;.
$$
\end{thm}

{\bf Remark : } For simplicity we are assuming here that $w_1$ and $w_2$ are H\"{o}lder continuous and bounded. 
These  assumptions could be relaxed by standard (but painfull) techniques. We do not know if the
uniform Lipschitz continuity assumption on $w_1$ with respect to $q$ and on $w_2$ with respect to $p$ can be relaxed.\\

As a consequence we have

\begin{thm} [Existence of a value] \label{ExistValue} Under assumptions (H), (\ref{Hypgij}) and (\ref{Isaacs}), the game has a value:
$$
V^+(t,x,p,q)=V^-(t,x,p,q)\qquad \forall (t,x,p,q)\in (0,T)\times \R^n\times\Delta(I)\times \Delta(J)\;.
$$
Furthermore $V^+=V^-$ is the unique solution in the dual sense of HJI equation (\ref{HJI}) with terminal condition
$$
V^+(T,x,p,q)=V^-(T,x,p,q)=\sum_{i=1}^I\sum_{j=1}^J p_iq_jg_{ij}(x)\qquad \forall (x,p,q)\in \R^n\times \Delta(I)\times \Delta(J)\;.
$$
\end{thm}

{\bf Proof of Theorem \ref{ExistValue} : } The Hamiltonian H defined by (\ref{Isaacs}) is known to satisfy (\ref{StructH}) (see \cite{CIL} for instance). 
>From the definition of $V^+$ and $V^-$ we have $V^-\leq V^+$. 
We have proved in Lemma \ref{reguV+V-} and Proposition  \ref{CvV+V-} that $V^+$ and $V^-$ are H\"{o}lder continuous, Lipschitz continuous 
with respect to $p$ and $q$, convex w.r. to $p$ and concave w.r. to $q$. From Corollary \ref{HJIV-} we know that $V^-$ is a supersolution
of (\ref{HJI}) in the dual sense while Proposition \ref{CasV+} states that $V^+$ is a supersolution of that same equation in the
dual sense. The comparison principle then states that $V^+\leq V^-$, whence the existence and the characterization  of the value:
$V^+=V^-$ is the unique solution in the dual sense of HJI equation (\ref{HJI}).
\cq

\subsection{Proof of the comparison principle}

The proof of Theorem \ref{comparison} relies on two arguments: first on a reformulation
of the notions of sub- and supersolutions by using sub- and superjets; second on
a new maximum principle described in the appendix.\\

Let us recall the notions of sub- and superjets of a function
$w: (0,T)\times \R^n\to \R$: the subjet $D^{2,-}w(\bart,\barx)$ is the set of 
$(\xi_t, \xi_x, X)\in \R^{n+1}\times {\cal S}_n$ such that
$$
w(t,x)\geq w(\bart, \barx)+\xi_t(t-\bart)
+\xi_x.(x-\barx)+\frac12X(x-\barx).(x-\barx)+o(|t-\bart|+|x-\barx|^2)\}
$$
and the superjet $D^{2,+}w$ is given by
$$
D^{2,+}w(\bart,\barx)=-D^{2,-}(-w)(\bart,\barx)
$$
When $w$ depends on other variables ($(p,q)$ or $(p, \hq)$ for instance), 
$D^{2,-}w$ and $D^{2,+}w$ always denote the sub- and superjets with respect to the $(t,x)$
variables only. For $w=w(t,x,p,\hq)$, we set
$$
\overline{D^{2,-}}w(\bart,\barx,p,\hq)=\left\{\begin{array}{c} (\xi_t, \xi_x, X)\in \R^{n+1}\times {\cal S}_n\;,\;
\exists (t_n,x_n, p_n,\hq_n)\to (\bart,\barx,p,\hat q),\\
\exists (\xi_t^n, \xi_x^n, X^n)\in D^{2,-}w(t_n,x_n,p_n,\hq_n)\\ {\rm with}\; 
 (\xi_t^n, \xi_x^n, X^n)\to (\xi_t, \xi_x, X)\end{array}\right\}\;.
$$
We use a symmetric notation for $\overline{D^{2,+}}w(\bart,\barx,\hp,q)$.

The following equivalent formulation of the notion of sub- and supersolution
is standard in viscosity solution theory, so we omit the proof:

\begin{prop}\label{DefAvecJets} A  map $w:(0,T)\times \R^n\times \Delta(I)\times \Delta(J)\to \R$
is a supersolution of equation (\ref{HJI}) in the dual sense if and only
if $w=w(t,x,p,q)$ is lower semicontinuous, concave with respect to $q$ and if, 
for any $(\bart,\barx,\hp,\barq)$ and any 
$(\xi_t, \xi_x, X)\in \overline{D^{2,+}}w^*(\bart,\barx, \hp,\barq)$ we have
$$
\xi_t-H(\bart, \barx, -\xi_x, -X, p, \barq)\geq 0
\qquad \forall p \in \partial^-_\hp w^*(\bart,\barx, \hp,\barq)\;.
$$
Symmetrically $w$ is a subsolution of (\ref{HJI}) in the dual sense if and only if
$w$ is upper semicontinuous, convex with respect to $p$ and if, for any 
$(\bart,\barx,\barp,\hq)$ and any $(\xi_t, \xi_x, X)\in \overline{D^{2,-}}w^\sharp(\bart,\barx, \barp,\hq)$ 
we have 
$$
\xi_t-H(\bart, \barx, -\xi_x, -X, \barp, q)\leq 0
\qquad \forall q \in \partial^+_\hq w^\sharp(\bart,\barx, \barp, \hq)\;.
$$
\end{prop}

\noindent {\bf Proof of Theorem \ref{comparison} :} 
Let us assume that
$$
\sup_{t,x,p,q}(w_1-w_2)>0\;.
$$
Since
$w_1$ and $w_2$ are H\"{o}lder continuous and bounded, classical arguments show that 
$$
M_{\epsilon, \eta, \alpha} :=\sup_{t,x,s,y,p,q} \left\{ w_1(t,x,p,q)-w_2(s,y,p,q)-
( \frac{|(t,x)-(s,y)|^2}{2\epsilon}+\frac{\alpha}{2}(|x|^2+|y|^2))+\eta t\right\}
$$
is finite and achieved at a point $(\bart, \barx, \bars, \bary, \barp_0, \barq_0)$.
One can also show that
\begin{equation}\label{Estimates0}
\lim_{\epsilon, \eta, \alpha\to0^+}M_{\epsilon, \eta, \alpha}=\sup_{t,x,p,q}(w_1-w_2)
\end{equation}
and that
\begin{equation}\label{Estimates1}
\frac{|(\bart,\barx)-(\bars,\bary)|^2}{\epsilon^2}\;,\; \alpha |\barx|^2\;,\; \alpha |\bary|^2 \leq 2M_\infty
\end{equation}
where $M_\infty=|w_1|_\infty+|w_2|_\infty$. Using (\ref{InitCond}) and the H\"{o}lder continuity of
$w_1$ and $w_2$ shows that $\bart<T$ and $\bars<T$ as soon as $\epsilon$, $\eta$ and $\alpha$ are small
enough.

>From the maximum principle (Theorem \ref{Ishii} stated in the Appendix), 
there are $(\barp, \barq)$, $(\hp, \hq)$ and  $X_1, X_2\in {\cal S}_n$ such that
$$
\barp \in \partial^-_\hp w_2^*(\bars, \bary,\hp, \barq), \;  \barq\in \partial^-_\hq w_2^\sharp (\bart, \barx,\barp, \hq)\;,
$$
$$
(-\frac{(\bart-\bars)}{\epsilon}+\eta, -  \frac{(\barx-\bary)}{\epsilon}-\alpha \barx, X_1)\in 
\overline{D^{2,-}}w_1^\sharp(\bars, \barx,\barp, \hq)\;, 
$$
$$
(\frac{(\bars-\bart)}{\epsilon},  \frac{(\bary-\barx)}{\epsilon}-\alpha\bary, X_2)\in \overline{D^{2,+}}w_2^*(\bars,\bary,\hp, \barq)
$$
and
\begin{equation}\label{Estimates2}
 \left(\begin{array}{ll} -X_1 & 0\\ 0 & X_2\end{array}\right) \leq (\frac{3}{\epsilon}+2\alpha)\left(\begin{array}{cc}I&-I\\-I& I\end{array}\right)
+(\alpha+\alpha^2\epsilon)I
\end{equation}
Since $w_1$ is a subsolution of (\ref{HJI}) in the dual sense and  $\barq\in \partial^-_\hq w_2^\sharp (\bart, \barx,\barp, \hq)$, Proposition \ref{DefAvecJets} states that
\begin{equation}\label{CondVisc1}
\eta-\frac{\bart-\bars}{\epsilon}-H\left(\bart, \barx, \frac{\barx-\bary}{\epsilon}+\alpha\barx,-X_1,\barp,\barq\right)\leq 0\;.
\end{equation}
In the same way, since
$w_2$ is a supersolution of (\ref{HJI}) in the dual sense and $\barp \in \partial^-_\hp w_2^*(\bars, \bary,\hp, \barq)$, we have
\begin{equation}\label{CondVisc2}
\frac{(\bars-\bart)}{\epsilon}-H\left(\bars, \bary, -\frac{(\bary-\barx)}{\epsilon}+\alpha\bary, -X_2, \barp,\barq\right)\geq 0\;,
\end{equation}

Using the structure condition (\ref{StructH}) on $H$, 
and plugging estimates (\ref{Estimates0}), (\ref{Estimates1}) and (\ref{Estimates2}) into (\ref{CondVisc1}) and (\ref{CondVisc2}) yields to a contradiction for $\epsilon$, $\alpha$ and $\eta$ sufficiently small  as in \cite{CIL}.
\cq

\section{Games with running cost}

We now investigate differential games with asymmetric information on the running cost and on the
terminal cost. The framework is basically the same as before. At the initial time, the cost
(now consisting in a running cost and a  terminal one) is chosen at random among
$I\times J$ possible costs. The index $i$ is announced to Player~I while the index $j$ is announced
to Player~II. Then the players play the game in order, for Player I to minimize the payoff and for Player II to maximize it.

In this section we keep the same terminology  and the same notations as in the previous part.
There is however a main difference: as we shall see later, 
in a game with a running cost, each player needs the knowledge of this running cost  to build his strategy. 
Since we assume that the running cost depends on the control of both players, this means that the players have to
observe the control of their opponent. This was not the case of the game before where the players only observed
the state of the system. For this reason we have to change the notion of strategy: 
in this section the notion of strategies introduced in Definition \ref{strat} is replaced by the following one:

\begin{Definition} 
A {\em strategy} for player I starting at time $t$ is a Borel-measurable map 
$\alpha : [t,T]\times\CR([t,T],\R^n)\times L^2([t,T],V)\rightarrow U$ for which
there exists $\delta>0$ such that, for all $s\in[t,T]$, $f,f'\in\CR([t,T],\R^n)$ and
$g,g'\in L^2([t,T],V)$, 
if $f=f'$ and $g=g'$ a.e. on $[t,s]$, then $\alpha(\cdot,f,g)=\alpha(\cdot,f',g')$ on  $[t,s+\delta]$.\\
We define strategies for player II in a symmetric way and denote by $\AR(t)$ (resp. 
$\BR(t)$) the set of strategies for player I (resp. player II).
\end{Definition}
We define random strategies as before (but with the modified notion of strategies) and still denote by $\AR_r(t)$ (resp. 
$\BR_r(t)$) the set of random strategies for player I (resp. player II). \\
We have an analogue of Lemma \ref{ptfix} :

\begin{Lemma}\label{ptfix2}
For all $(t,x)$ in $[0,T]\times\R^n$, for all $(\alpha,\beta)\in\AR(t)\times\BR(t)$, 
there exists a unique couple of controls $(u,v)\in\UR(t)\times\VR(t)$ that satisfies $P-$a.s.
\begin{equation}
 \label{alphau2}
(u,v)=(\alpha(\cdot,X^{t,x,u,v}_\cdot,v_\cdot),\beta(\cdot,X^{t,x,u,v}_\cdot,u_\cdot)) 
\mbox{ a.e. on }
[t,T].
\end{equation}

\end{Lemma}

One can easily check that the results of the previous parts (i.e., \ref{reguV+V-}, Proposition  \ref{CvV+V-}, Corollary \ref{HJIV-}  and  Proposition \ref{CasV+})
still hold true with the modified notion of strategy. In particular, the game with terminal payoff studied before has a value.\\
Let us fix $I,J\in \N$. For $1\leq i\leq I$ and $1\leq j\leq J$ we consider the terminal cost $g_{ij}:\R^n\to \R$ and 
 the running cost $\ell_{ij}:[0,T]\times \R^n\times U\times V\to \R$ on which we do the following assumptions:
\begin{equation}\label{Hyplijgij}
\begin{array}{c}
\mbox{\rm For any $1\leq i\leq I$ and $1\leq j\leq J$, 
$\ell_{ij}$ and $g_{ij}$ are continuous in all variables,}\\ 
\mbox{\rm uniformly Lipschitz continuous with respect to $x$ and bounded.}
\end{array}
\end{equation}
For fixed $(i,j)\in\{ 1,\ldots,I\}\times\{ 1,\ldots,J\}$ and strategies 
$(\alpha,\beta)\in\AR(t)\times\BR(t)$, we set 
\[ J_{ij}(t,x,\alpha,\beta)=E[\int_t^T\ell_{ij}(s,X^{t,x,\alpha,\beta}_s, \alpha_s, \beta_s)ds+ g_{ij}(X^{t,x,\alpha,\beta}_T)]\;,\]
where as before $(\alpha, \beta)$ denotes the unique pair of controls such that (\ref{alphau2}) holds. \\
The payoff of two random strategies $(\overline\alpha,\overline\beta)\in \AR_r(t)\times\BR_r(t)$, with
$\overline\alpha=(\alpha^1,\ldots,\alpha^R;r^1,\ldots,r^R)$ and
$\overline\beta=(\beta^1,\ldots,\beta^S;s^1,\ldots,s^S)$, 
 is the average of the payoffs with respect to the probability distributions 
associated to the strategies:
\[J_{ij}(t,x,\overline\alpha,\overline\beta)=
\sum_{k=1}^R\sum_{l=1}^Sr^ks^lE[\int_t^T\ell_{ij}(s,X^{t,x,\alpha^k,\beta^l}_s, \alpha^k_s, \beta_s^l)ds+ g_{ij}(X^{t,x,\alpha^k,\beta^l}_T)].\]
Finally, the payoff of the game is, for 
$(\hat\alpha,\hat\beta)=((\bar \alpha_i)_{1\leq i\leq I}, (\bar \beta_j)_{1\leq j\leq J}) \in(\AR_r(t))^I\times(\BR_r(t))^J$,
\[ J^{p,q}(t,x,\hat\alpha,\hat\beta)=\sum_{i=1}^I\sum_{j=1}^Jp_iq_j
J_{ij}(t,x,\overline\alpha_i,\overline\beta_j).\]
We define the value functions for the game with running cost as before by
\[\begin{array}{c}
V^+(t,x,p,q)=\inf_{\hat\alpha\in(\AR_r(t))^I}\sup_{\hat\beta\in
\BR_r(t))^J}J^{p,q}(t,x,\hat\alpha,\hat\beta),\\
V^-(t,x,p,q)=\sup_{\hat\beta\in
(\BR_r(t))^J}\inf_{\hat\alpha\in(\AR_r(t))^I}J^{p,q}(t,x,\hat\alpha,\hat\beta).
\end{array}\]

In our game with running cost, Isaacs' assumption takes the following form: for all  $(t,x)\in[0,T]\times
\R^n$, $(p, q)\in \Delta(I)\times \Delta(J)$, $\xi \in\R^n$,  and all $A\in {\cal S}_n$:
\begin{equation}\label{IsaRun}
\begin{array}{l}
\inf_u\sup_v\{ <b(t,x,u,v),\xi>+\frac12
Tr(A\sigma(t,x,u,v)\sigma^*(t,x,u,v))-\sum_{i,j} \ell_{ij}(t,x,u,v)p_iq_j\}=\\
\hspace{1cm}\sup_v\inf_u\{ <b(t,x,u,v),\xi>+\frac12
Tr(A\sigma(t,x,u,v)\sigma^*(t,x,u,v))-\sum_{i,j} \ell_{ij}(t,x,u,v)p_iq_j\}
\end{array}
\end{equation}
We set 
\[\begin{array}{rl}
H(t,x,\xi,A,p,q)=&\inf_u\sup_v\{ <b(t,x,u,v),\xi>\\&+\frac12
Tr(A\sigma(t,x,u,v)\sigma^*(t,x,u,v))-\sum_{i,j} \ell_{ij}(t,x,u,v)p_iq_j\}\;.
\end{array}\]

\begin{thm} \label{ValueRun} Assume that (H), (\ref{Hyplijgij}) and (\ref{IsaRun}) hold. 
Then the game has a value: $V^+=V^-$, which is the unique solution in the dual sense
of 
\begin{equation}\label{HIJrun}
w_t+H(t,x,Dw,D^2w,p,q)=0
\end{equation}
with terminal condition
$$
V^+(T,x,p,q)=V^-(T,x,p,q)=\sum_{i=1}^I\sum_{j=1}^J p_iq_jg_{ij}(x)\qquad \forall (x,p,q)\in \R^n\times \Delta(I)\times \Delta(J)\;.
$$
\end{thm}

In order to prove Theorem \ref{ValueRun} it will be convenient to have the following equivalent definition of dual solutions
of the Hamilton-Jacobi equation (\ref{HIJrun}):

\begin{Lemma}\label{CaracSuperSol} Let $w:[0,T]\times \R^N\times \Delta(I)\times \Delta(J)\mapsto \R$ be lower-semicontinuous, 
uniformly Lipschitz continuous with respect to $p$ and concave 
with respect to $q$. Then the following statements are equivalent:
\begin{itemize}
\item[(i)] $w$ is a dual supersolution of (\ref{HIJrun}).

\item[(ii)] for any $(\hat p, \bar q)\in \R^I\times \Delta(J)$, for any ${\cal C}^2$ test function $\phi=\phi(t,x)$
such that 
$$
(t,x,p)\mapsto w(t,x,p,\bar q)-\phi(t,x)-\langle\hat p,p\rangle
$$
has a global minimum at some point $(\bar t,\bar x,\bar p)\in [0,T)\times \R^N\times \Delta(I)$, we have
\begin{equation}\label{ConSurSol}
\phi_t(\bar t,\bar x)+H(\bar t,\bar x,D\phi(\bar t,\bar x,\bar p),D^2\phi(\bar t,\bar x,\bar p), \bar p,\bar q) \leq 0\;.
\end{equation}
\end{itemize}
\end{Lemma}

\noindent{\bf Remark : } A symmetric statement holds for dual subsolutions.\\

{\bf Proof of Lemma \ref{CaracSuperSol} : } Let us assume that $w$ is a supersolution 
and let $\phi\in {\cal C}^2$, $(\hat p,\bar q)\in \R^I\times\Delta(J)$ such that 
\begin{equation}\label{unmin}
(t,x,p)\mapsto w(t,x,p,\bar q)-\phi(t,x)-\langle\hat p,p\rangle
\end{equation}
 has a global minimum at some point $(\bar t,\bar x,\bar p)\in [0,T)\times \R^N\times \Delta(I)$.
We note that this implies that $\bar p\in \partial^-_{\hat p}w^*(\bar t, \bar x,\hat p, \bar q)$. Moreover,
 taking the supremum over $p$ in (\ref{unmin}), we have that $(t,x)\to -\phi(t,x)-w^*(t,x,\hat p,\bar q)$ has  a global minimum at $(\bar t,\bar x)$.
Since $w^*$ is a subsolution of the dual equation, we get
$$
-\phi_t+H^*(\bar t,\bar x, -D\phi, -D^2\phi, \bar p, \bar q)\geq 0
$$
at $(\bar t, \bar x)$, because $\bar p\in \partial^-_{\hat p}w^*(\bar t, \bar x,\hat p, \bar q)$. Whence inequality (\ref{ConSurSol}).

Conversely let us assume that $w$ satisfies (ii). Let $\phi$ be a ${\cal C}^2$ test function such that $(t,x)\to w^*(t,x,\hat p,\bar q)-\phi(t,x)$
has a maximum at some point $(\bar t, \bar x)\in (0,T)\times \R^N$ for some $(\hat p,\bar q)\in \R^I\times \Delta(J)$.
Without loss of generality we can assume that this maximum is a global one.  Let $\bar p\in  \partial^-_{\hat p}w^*(\bar t, \bar x,\hat p, \bar q)$.
>From the definition of $w^*$, we also have that 
$(t,x,p)\to \langle\hat p,p\rangle - w(t,x,p, \bar q)-\phi(t,x)$ has a global maximum at $(\bar t, \bar x, \bar p)$, i.e.,  
$(t,x,p)\to w(t,x,p, \bar q)+\phi(t,x)-\langle\hat p,p\rangle $ has a global minimum at $(\bar t, \bar x, \bar p)$. From (\ref{ConSurSol}) we get
$$
-\phi_t+H(\bar t,\bar x, -D\phi, -D^2\phi, \bar p, \bar q)\leq 0\;,
$$
the desired inequality.
\cq

{\bf Proof of Theorem \ref{ValueRun} : } Following standard arguments, one first checks that 
$V^+$ and $V^-$ are globally H\"{o}lder continuous, and uniformly Lipschitz continuous
with respect to $p$ and $q$. In order to prove other properties of the value functions,
let us introduce an extended differential game in $\R^{n+IJ}$. This game with asymmetric information and terminal payoff is
defined by the dynamics
\begin{equation}\label{dynrun}
\begin{array}{l}
dX_s=b(s,X_s,u_s,v_s)ds+\sigma(s,X_s,u_s,v_s)dB_s,\; s\in[t,T],\\
dZ_{ij,s}=\ell_{ij}(s,X_s,u_s,v_s)ds,\\
X_t=x,\; Z_{ij,t}=z_{ij},
\end{array}
\end{equation}
where $(t,x,z)\in[0,T]\times\R^n\times \R^{IJ}$, with $z=(z_{ij})$, and the terminal 
$\tilde g_{ij}(x,z)=z^{ij}+g_{ij}(x)$. We denote by $\tilde V^+$ and $\tilde V^-$ the upper and lower
value of this game. We note that 
\begin{equation}\label{VTildeV}
\tilde V^\pm(t,x,z,p,q)=V^\pm(t,x,p,q)+\sum_{ij}z_{ij}p_iq_j\;.
\end{equation}

Following the proofs of Proposition  \ref{CvV+V-}, one can check that
$\tilde{V}^+$ and $\tilde{V}^-$ are convex in $p$ and concave in $q$. Hence so are $V^+$ and $V^-$.
As in Corollary \ref{HJIV-}  and  Proposition \ref{CasV+}, one can also show that $\tilde{V}^{-}$ is a dual supersolution of the HJ equation
$$
\tilde{w}_t+\tilde{H}^-(t,x,z,D_{x,z}w, D^2_{x}w)=0
$$
where, for $(t,x,z)\in \R^{n+IJ}$, $\xi_x\in \R^n$, $\xi_z\in \R^{IJ}$ and $A\in {\cal S}_n$,
\[\begin{array}{rl}
\tilde H^-(t,x, z, \xi_x, \xi_z, A)=& \sup_{v\in V}\inf_{u\in U}\{ <b(t,x,u,v),\xi_x>\\
&+\frac12
Tr(A\sigma(t,x,u,v)\sigma^*(t,x,u,v))+\sum_{i,j} \ell_{ij}(t,x,u,v)\xi_{z,ij}\},
\end{array}\]
while $\tilde{V}^{+}$ is a dual subsolution of the HJ equation
$$
\tilde{w}_t+\tilde{H}^+(t,x,z,D_{x,z}w, D^2_{x}w)=0
$$
where 
\[\begin{array}{rl}
\tilde H^+(t,x, z, \xi_x, \xi_z, A)=& \inf_{u\in U}\sup_{v\in V}\{ <b(t,x,u,v),\xi_x>\\
&+\frac12
Tr(A\sigma(t,x,u,v)\sigma^*(t,x,u,v))+\sum_{i,j} \ell_{ij}(t,x,u,v)\xi_{z,ij}\}.
\end{array}\]
Note that this is precisely at this point that the players have to use the new definition of strategies.
Indeed, in order to build their strategies in the sub- and superdynamic programming, they have to compute the running costs $Z_{ij}$
(see the proof of Proposition \ref{dynprog}).
This is possible since, at time $s$, they know the controls $u_\cdot$ and $v_\cdot$ and the trajectory $X_\cdot$ up to time $s-\delta$, and therefore can
compute $Z_{ij,s}=z_{ij}+\int_t^s \ell_{ij}(\tau, X_\tau, u_\tau, v_\tau)d\tau$.

Using Lemma \ref{CaracSuperSol} one can then show that $V^{-}$ is a dual supersolution of the HJ equation
$$
w_t+H^-(t,x,Dw, D^2w, p,q)=0
$$
where 
\[\begin{array}{rl}
H^-(t,x,  \xi, A, p,q)=& \sup_{v\in V}\inf_{u\in U}\{ <b(t,x,u,v),\xi_x>\\
&+\frac12
Tr(A\sigma(t,x,u,v)\sigma^*(t,x,u,v))+\sum_{i,j} \ell_{ij}(t,x,u,v)p_iq_j\}
\end{array}\]
while $V^+$ is a dual subsolution of the HJ equation
$$
w_t+H^+(t,x,Dw, D^2w, p,q)=0
$$
where 
\[\begin{array}{rl}
H^+(t,x,  \xi, A, p,q)=& \inf_{u\in U}\sup_{v\in V}\{ <b(t,x,u,v),\xi_x>\\
&+\frac12
Tr(A\sigma(t,x,u,v)\sigma^*(t,x,u,v))+\sum_{i,j} \ell_{ij}(t,x,u,v)p_iq_j\}
\end{array}\]
Finally combining Isaacs' assumption, which states that $H:=H^+=H^-$, the fact that $H$ satisfies assumption
(\ref{StructH}) and the comparison principle shows that $V^+=V^-$
is the unique dual solution of (\ref{HIJrun}).
\cq

\section{Appendix : A maximum principle}

The following result---used in a crucial way in the proof of the comparison principle---is 
an adaptation to our framework of the maximum principle for semicontinuous functions (see Theorem 3.2 of \cite{CIL}):

\begin{thm}[Maximum principle]\label{Ishii} For $k=1,2$, let ${\cal O}_k$ be open subsets of $\R^{n_k}$ and $w_k:{\cal O}_k\times \Delta(I)\times \Delta(J)\to \R$
be such that
\begin{itemize}
\item[(i)] $w_1=w_1(x,p,q)$ is upper semicontinuous in all variables,  convex with respect to $p$ and uniformly Lipschitz continuous with respect to $q$, 
\item[(ii)] $w_2=w_2(y,p,q)$ is lower semicontinuous in all its variables, concave with respect to $q$
and uniformly Lipschitz continuous with respect to $p$,
\item[(iii)] there is some ${\cal C}^2$ map $\phi: {\cal O}_1\times {\cal O}_2\to \R$ and some point 
$(\barx,\bary)\in {\cal O}_1\times {\cal O}_2$ such that the map
$$
(x,y)\to \max_{p,q}\left\{w_1(x,p,q)-w_2(y,p,q) -\phi(x,y)\right\}
$$
has a maximum at $(\barx,\bary)$. 
\end{itemize}

Then, for any $\epsilon>0$, there are $(\barp,\barq) \in \Delta(I)\times \Delta(J)$, $(\hp,\hq)\in\R^I\times \R^J$ and $(X_1, X_2)\in {\cal S}_{n_1}\times {\cal S}_{n_2}$ such that
the map
$$
(x,y,p,q)\to  w_1(x,p,q)-w_2(y,p,q) -\phi(x,y)
$$
has a maximum at $(\barx,\bary,\barp,\barq)$, 
\begin{equation}\label{Ishii0}
\barp \in \partial^-_\hp w_2^*(\bary,\hp, \barq), \;  \barq\in \partial^-_\hq w_1^\sharp (\barx,\barp, \hq)\;,
\end{equation}
\begin{equation}\label{Ishii1}
(-D_{x}\phi(\barx, \bary) , X_1)\in \overline{D^{2,-}}w_1^\sharp(\barx,\barp, \hq), \; (D_{y}\phi(\barx, \bary) , X_2)\in \overline{D^{2,+}}w_2^*(\bary,\hp, \barq)
\end{equation}
and
\begin{equation}\label{Ishii2}
\left(\frac{1}{\epsilon}+\|A\|\right) I \leq \left(\begin{array}{ll} -X_1 & 0\\ 0 & X_2\end{array}\right) \leq A+\epsilon A^2
\end{equation}
with $A= D^2\phi(\barx, \bary)$.
\end{thm}

{\bf Remark : } Compared with the classical maximum principle, the additional difficulty here is the
fact that we need elements of $\overline{D^{2,-}}w_1^\sharp$ and of $\overline{D^{2,+}}w_2^*$
while we have only information on the behavior of the difference $w_1-w_2-\phi$.\\

{\bf Proof of Theorem \ref{Ishii} : } We follow closely the proof of Theorem 3.2 of \cite{CIL}. Let us start by some reductions:

\noindent {\bf Reductions : } As in \cite{CIL}, we can assume without loss of generality 
that ${\cal O}_k=\R^{n_k}$, $\barx=\bary=0$ and $\phi(x,y)=A(x,y).(x,y)$ and 
\begin{equation}\label{max0}
\max_{x,y,p,q} \{ w_1(x,p,q)-w_2(y,p,q)-\phi(x,y)\}=0\;. 
\end{equation}
We can also assume that, for any $(\barp', \barq')\in \Delta(I)\times \Delta(J)$, 
\begin{equation}\label{HypInterieur}
\begin{array}{c}
\mbox{\rm if $(\barx,\bary, \barp', \barq')$ is a maximum point of $w_1-w_2-\phi$,}\\
\mbox{\rm  then $(\barp', \barq')$ belongs to the interior of $\Delta(I)\times \Delta(J)$.}
\end{array}
\end{equation}
Indeed, let us assume that Theorem \ref{Ishii} holds true under this additionnal assumption
and let us prove that it holds true without.
Let
$$
z(x,y,p,q)= w_1(x,p,q)-w_2(y,p,q) -\phi(x,y)\;.
$$
Among the $(p,q)$ for which $z(\barx,\bary,p,q)$ has a maximum, let us choose
$(\barp_0,\barq_0)$ such that the total number of indices $i$ and $j$ for which
$(\barp_0)_i=0$ or $(\barq_0)_j=0$ is maximal.
Let us denote by $I'$ and $J'$  the set of indices $i$ and $j$ 
for which $(\barp_0)_i>0$ and $(\barq_0)_j>0$. We then define
$w_1'$, $w_2'$, $z'$, $\barp_0'$ and $\barq_0'$ as the natural restriction of $w_1$, $w_2$, $z$, 
$\barp_0$ and $\barq_0$ to $\Delta(I')$ and $\Delta(J')$. 
We note that $(\barx,\bary,\barp_0',\barq_0')$ is a maximum point of $z$ on $\R^{n_1+n_2}\times
\Delta(I')\times \Delta(J')$ and that assumption (\ref{HypInterieur})
holds, since otherwise one would have a contradiction with the particular
choice of $(\barp_0,\barq_0)$. 

Using now Theorem \ref{Ishii} with assumption (\ref{HypInterieur}), we can build 
$(\barp',\barq') \in \Delta(I')\times \Delta(J')$, $\hq'\in \partial^+_q w_1(\barx,\barp', \barq')$,  
$\hp'\in \partial^-_p w_2(\bary,\barp', \barq')$ and $(X_1, X_2)\in {\cal S}_{n_1}\times {\cal S}_{n_2}$ such that 
(\ref{Ishii0}), (\ref{Ishii1}) and (\ref{Ishii2}) hold. Then we extend  $(\barp',\barq')$ to $ (\barp,\barq) \in \Delta(I)\times \Delta(J)$ by setting $\barp_i=\barp_i'$ for $i\in I'$ and $\barp_i=0$ otherwise,  
and $\barq_j=\barq_j$ for $j\in J'$ and $\barq_j=0$
otherwise. We also extend $\hq'$ to $\hq\in \partial^+_q w_1(\barx,\barp, \barq)$
and $\hq'$ to $\hq\in \partial^+_q w_1(\barx,\barp, \barq)$ by setting $\hq_j=M$ for $j\in J\backslash J'$ and $\hp_i=-M$ for $i\in I\backslash I'$, 
where $M$ is a Lipschitz constant of $w_1$ and $w_2$ with respect to $q$ and $p$ respectively. This defines
$\barp,\barq, \hq, \hp$ and $(X_1, X_2)$ for which (\ref{Ishii0}), (\ref{Ishii1}) and (\ref{Ishii2}) hold.

So it remains to prove that Theorem \ref{Ishii} holds true under the additional assumption (\ref{HypInterieur}).\\

\noindent {\bf Step 1 : introduction of the inf- and supconvolutions. } As in \cite{CIL}, we have
$$
(w_1(x',p,q)-\frac{\lambda}{2}|x'-x|^2)-(w_2(y',p,q)-\frac{\lambda}{2}|y'-y|^2) \leq \langle(A+\epsilon A^2)(x,y),(x,y)\rangle
$$
for any $(x,x',y,y',p,q)$, where $\lambda= \frac{1}{\epsilon}+\|A\|$. Let us set for $\lambda'\in (0, \lambda)$, 
$$
\hw_1(x,p,q)=\max_{x'\in \R^{n_1},\;q'\in \Delta(J)}(w_1(x',p,q')-\frac{\lambda}{2}|x'-x|^2-\frac{\lambda'}{2}|q'-q|^2)
$$
and
$$
\hw_2(y,p,q)=\min_{y'\in \R^{n_2},\; p'\in \Delta(I)}(w_2(y',p',q)+\frac{\lambda}{2}|y'-y|^2+\frac{\lambda'}{2}|p'-p|^2)\;.
$$
With these definition we have that $\hw_1$ is semiconvex in all its variables 
with a modulus $\lambda'$,  semiconvex in $x$ with a modulus $\lambda$ and convex in $p$ (because $w_1$ is convex in $p$
by assumption). In the same way,  $\hw_2$ is semiconcave in all its variables with a modulus $\lambda'$, semiconvex in $y$ with a modulus $\lambda$ 
and concave in $q$ (because $w_2$ is concave in $q$
by assumption). Moreover
\begin{equation}\label{maxhw}
\hw_1(x,p,q)-\hw_2(y,p,q)- \langle (A+\epsilon A^2)(x,y),(x,y)\rangle\leq 0\qquad \forall (x,y,p,q)\;.
\end{equation}
Since $w_1\leq \hw_1$ and $w_2\geq \hw_2$, there are some $(p,q)$ such that equality holds
in (\ref{maxhw}) at $(0,0,p,q)$. Furthermore, if equality holds at $(0,0,p,q)$, then $(0,0,p,q)$ is a maximum point in (\ref{max0}) and assumption (\ref{HypInterieur}) states that
$(p,q)$ belongs to the interior of $\Delta(I)\times \Delta(J)$.\\

\noindent {\bf Step 2 : use of Jensen maximum principle.} Let us now introduce some small pertubation of the equation: for $\alpha>0$
and $\zeta=(\zeta_x,\zeta_y, \zeta_p, \zeta_q)\in \R^{n_1+n_2+I+J}$, we set
$$
\begin{array}{rl}
z_\zeta(x,y,p,q)\; =\ & \hw_1(x,p,q)-\hw_2(y,p,q)- (A+\epsilon A^2)(x,y).(x,y)\\
& \qquad -\alpha(|x|^2+|y|^2+|p|^2-|q|^2) -\langle \zeta,(x,y,p,q)\rangle\;.
\end{array}
$$
Note that, because of the penalisation term $\alpha(|x|^2+|y|^2)$, for any $\eta>0$, we can choose $\gamma$
small enough such that, for any $\zeta$ such that $|\zeta|\leq \gamma$, any maximum of $z_\zeta$ is of the form $(x,y,p,q)$ for some $(x,y)\in B_\eta$.

Let $\gamma$ as above. 
Since $z_0$ is semiconvex, has a maximum at $(0,0, p,q)$, Jensen maximum principle (see Lemma A.3 of \cite{CIL} for instance) states that the set
$$
 E_\gamma =\left\{\begin{array}{l}  (x,y,p,q)\in B_\eta\times \Delta(I)\times \Delta(J) \;,\;
 \exists \zeta \;,\;  |\zeta|\leq \gamma\;,\;  \mbox{\rm such that }\\ 
\qquad (i) \quad z_{\zeta}\;  \mbox{\rm has a maximum at $(x,y,p,q)$ and }\\
\qquad (ii) \quad \mbox{\rm  $\hw_1$ and $\hw_2$ have a derivative at $(x,y,p,q)$}\end{array} \right\}
$$
has a positive measure. We note that in the quoted Lemma A.3, the maximum is required to be strict ;
this assumption is only used in \cite{CIL} to localize the maximum points, which is not needed here.\\

We also note for later use that, if $(x,y, p,q)\in E_\gamma$, there is some $\zeta=(\zeta_x,\zeta_y,\zeta_p,\zeta_q)$ with $|\zeta|\leq \gamma$ such that
$z_\zeta$ has a maximum at $(x,y, p,q)$. In particular, this implies that 
$$
q'\to \hw_1(x,p,q')-\hw_2(y,p,q')+\alpha |q'|^2-\langle \zeta_q, q'\rangle
$$
has a maximum at $q$. Since $\hw_2$ is concave in $q$, $\hw_1$ coincides with its concave hull
with respect to $q$ at $(x,p,q)$. Hence, if we set $\hq=\frac{\partial \hw_1(x,p,q)}{\partial q}$, then
\begin{equation} \label{Pptehq}
\hw_1(x,p,q)+\hw_1^\sharp(x,p,\hq)=q.\hq\; {\rm  and }\; 
q\in \partial^+_\hq\hw_1^\sharp(x,p,\hq)\;.
\end{equation}
In the same way, if  we set $\hp=\frac{\partial \hw_2(y,p,q)}{\partial p}$, then we have
\begin{equation} \label{Pptehp}
\hw_2(x,p,q)+\hw_2^*(y,\hp,q)=p.\hp\; {\rm  and }\; 
p\in \partial^-_\hp\hw_2^*(x,\hp,q)\;.
\end{equation}

\noindent {\bf Step 3 : measure estimate of a subset of $E_\gamma$.} 
Let $E'_\gamma$ be the set of points $(x,y,p,q)\in E_\gamma$ such that
$\hw_1^\sharp$ has a second order Taylor expansion at $(x,p, \frac{\partial \hw_1}{\partial q}(x, p,q))$ and 
$\hw_2^*$ has a second order Taylor expansion at 
$(y, \frac{\partial \hw_2}{\partial p}(x, p,q),q)$. 
Our aim is to show that $E'_\gamma$ has a full measure in $E_\gamma$.\\

For this we note that $E'_\gamma=E^1_\gamma\cap E^2_\gamma$ where
$$
E^1_\gamma=\left\{\begin{array}{l} (x,y,p,q)\in E_\gamma\;,\;  \hw_1^\sharp\; \mbox{\rm has a second order Taylor expansion}\\
\qquad  \mbox{\rm  at $(x,p, \frac{\partial \hw_1}{\partial q}(x, p,q))$}\end{array}\right\}
$$
and
$$
E^2_\gamma=\left\{\begin{array}{l} (x,y,p,q)\in E_\gamma\;,\;  \mbox{\rm $\hw_2^*$ has a second order Taylor expansion}\\
\qquad \mbox{\rm  at $(y, \frac{\partial \hw_2}{\partial p}(x, p,q),q)$}\end{array} \right\}
$$
It is therefore enough to show that $E^1_\gamma$ and $E^2_\gamma$ have a full measure in $E_\gamma$.
We only do the proof for $E^1_\gamma$, the proof for $E^2_\gamma$ being symmetric.

Let us set, for any $(x,y,p)$,  
$$
E_\gamma(x,y,p)=\{ q \in \Delta(J)\;,\; (x,y,p,q)\in E_\gamma \}
$$
and
$$
E^1_\gamma(x,y,p)=\{ q \in \Delta(J)\;,\; (x,y,p,q)\in E^1_\gamma \}
$$
Since $E_\gamma$ has a positive measure, from Fubini Theorem we have to show that, for any $(x,y,p)$ such that the set
$E_\gamma(x,y,p)$ has a positive measure, the set $E^1_\gamma(x,y,p)$
has a full measure in $E_\gamma(x,y,p)$. \\

For this, let us introduce the map $\Phi: q\to \frac{\partial \hw_1(x,p,q)}{\partial q}$ defined on $E_\gamma(x,y,p)$. We are going to show that
\begin{equation}\label{InegMesure}
\forall q_1,q_2\in E_\gamma(x,y,p), \; |q_1-q_2|\leq \frac{1}{2\alpha} |\Phi(q_1)-\Phi(q_2)|\;,
\end{equation}
which will imply that
\begin{equation}\label{Mesure}
\forall E\subset E_\gamma(x,y,p)\; {\rm measurable},\; {\cal L}^J(E)\leq \frac{1}{(2\alpha)^I} {\cal L}^J(\Phi(E))\;,
\end{equation}
where ${\cal L}^J$ denotes the Lebesgue measure in $\R^J$.
Then we will prove that (\ref{Mesure}) implies our claim.\\

{\it Proof of (\ref{InegMesure}) : } Let $q_1,q_2\in E_\gamma(x,y,p)$. There are $\zeta_1$ and $\zeta_2$ such that 
$z_{\zeta_k}$ has a maximum at $(x,y,p,q_k)$ for $k=1,2$. The first order optimality conditions
imply that
$$
\Phi(q_k)=\frac{\partial \hw_2(y,p,q_k)}{\partial q}-2\alpha q_k+\zeta_{k,q}\qquad {\rm for}\; k=1,2.
$$
Using again the optimality of $z_{\zeta_1}$ at $q_1$ and the fact that $q\to \hw_2(y,p,q)$ is concave, we have
$$
\begin{array}{rl}
\hw_1(x,p,q_2)\; \leq & \hw_1(x,p,q_1)+\langle \left(\frac{\partial \hw_2(y,p,q_k)}{\partial q}-2\alpha q_1+\zeta_{1,q}\right),(q_2-q_1)\rangle-\alpha |q_2-q_1|^2\\
\leq & \hw_1(x,p,q_1)+\langle \Phi(q_1),(q_2-q_1)\rangle-\alpha |q_2-q_1|^2
\end{array}
$$
Reversing the role of $q_1$ and $q_2$ gives
$$
\hw_1(x,p,q_1)\leq \hw_1(x,p,q_2)+\langle \Phi(q_2),(q_1-q_2)\rangle-\alpha |q_2-q_1|^2
$$
Adding the two previous inequalities then leads to
$$
0\leq (\Phi(q_2)-\Phi(q_1)).(q_1-q_2)-2\alpha|q_2-q_1|^2\;.
$$
Whence (\ref{InegMesure}).\\

{\it Proof of (\ref{Mesure}) : } Let $E$ be a measurable subset of $E_\gamma(x,y,p)$. We note that (\ref{InegMesure}) states that
$\Phi$ is a bijection between $E$ and its image, with a $\frac{1}{2\alpha}-$Lipschitz continuous inverse. Hence 
$$
{\cal L}^I(E)={\cal L}^I(\Phi^{-1}(\Phi(E)))\leq \frac{1}{(2\alpha)^I} {\cal L}^I(\Phi(E))\;,
$$
i.e., (\ref{Mesure}) holds. \\

We finally show that $E^1_\gamma(x,y,p)$ has a full measure in $E_\gamma(x,y,p)$ for any $(x,y,p)$ such that
$E_\gamma(x,y,p)$ has a positive measure. Let $F$ be the set of $(x,p,\hq)$ such that
$\hw_1^\sharp$ has a second order Taylor expansion at $(x,p,\hq)$. Since $F$ has a full measure, for almost all
$(x,p)\in \R^n\times \Delta(I)$, the set $F(x,p)=\{\hq\in \R^J\;,\; (x,p,\hq)\in F\}$ has a full measure 
in $\R^J$. Let $(x,p)$ be such a pair and such that $E_\gamma(x,y,p)$ has a positive measure. Then 
$\Phi(E_\gamma(x,y,p))$ also has a positive measure from (\ref{Mesure}). Since $\Phi(E_\gamma(x,y,p))\backslash F(x,p)$
has a zero measure and since
$$
\Phi^{-1}\left(\Phi(E_\gamma(x,y,p))\backslash F(x,p)\right)= E_\gamma(x,y,p)\backslash E^1_\gamma(x,y,p)\;,
$$
using again (\ref{Mesure}) shows that $ E_\gamma(x,y,p)\backslash E^1_\gamma(x,y,p)$ has a zero measure. This completes our claim.\\

\noindent {\bf Step 4 : (further) magic properties of sup-convolution. } We now explain 
that one can use second order Taylor expansions of $\hw_1^\sharp$ and $\hw_2^*$ to get elements of 
$D^{2,-}w_1^\sharp$, $D^{2,+}\hw_1$, $D^{2,+}w_2^*$ and $D^{2,-}\hw_1$.

>From our assumption (\ref{HypInterieur}), we know that, for $\epsilon$ small enough,
if $(0,0, p,q)$ realizes the equality in (\ref{maxhw}), then $(p,q)$ belongs to the interior of $\Delta(I)\times \Delta(J)$. 
Hence we can find $\alpha, \gamma>0$ so small that,
for any $\zeta$ with $|\zeta|\leq \gamma$, if $(x,y,p,q)$ realizes the maximum of $z_\zeta$, then 
$(p,q)$ belongs to the interior of $\Delta(I)\times \Delta(J)$.

Let us now fix $\gamma>0$ small enough and let us compute $\hw_1^\sharp$ at $(x,p,\hq)$ for $(x,y, p,q)\in E'_\gamma$ and $\hq=\frac{\partial \hw_1(x,p,q)}{\partial q}$. We have
\begin{equation}\label{hw1sharp}
\hw_1^\sharp(x,p,\hq)\; = \;  \min_{x',q',q"} (q'.\hq+\frac{\lambda}{2}|x'-x|^2+\frac{\lambda'}{2}|q"-q'|^2)-w_1(x',p,q")) 
\end{equation}
>From (\ref{Pptehq}), we have that  $\hw_1(x,p,q)+\hw_1^\sharp(x,p,\hq)=q.\hq$ and $q\in \partial^+_\hq\hw_1^\sharp(x,p,\hq)$. 
In particular, $q'=q$ is a minimum
point in (\ref{hw1sharp}). Since $q$ belongs to the interior of $\Delta(J)$, the optimality conditions
 imply that, if $(x',q,q")$ is a minimum of (\ref{hw1sharp}), then  $q=q"-\frac{1}{\lambda'}\hq$. Therefore
$$
\begin{array}{rl}
\hw_1^\sharp(x,p,\hq)\; =& -\frac{1}{2\lambda'}|\hq|^2+\min_{x',q"} (q".\hq-w_1(x',p,q")+\frac{\lambda}{2} |x'-x|^2) \\
= & -\frac{1}{2\lambda'}|\hq|^2+\min_{x'} (w_1^\sharp(x',p,\hq)+\frac{\lambda}{2} |x'-x|^2) \;,
\end{array}
$$
In particular, $q"\in \partial^+_\hq w_1(x', p,\hq)$, which shows that 
\begin{equation}\label{charcq}
q+\frac{1}{\lambda'}\hq\in \partial^+_\hq w_1^\sharp(x+\xi/\lambda, p,\hq)
\end{equation}
Moreover, $x\to \hw_1^\sharp(x,p,\hq)$ is equal, up to a constant, to the inf-convolution of $w_1^\sharp$ with respect to $x$. Since $\hw_1^\sharp$ has a 
second order Taylor expansion in $x$ at such a point $(x,p,\hq)$, the classical ``magic properties" of inf-convolution
(see Lemma A.4 of \cite{CIL}) state that $x'=x+\xi/\lambda$ and 
\begin{equation}\label{D2-w1sharp}
(D\hw^\sharp_1(x,p,\hq), D^2\hw^\sharp_1(x,p,\hq))\in D^{2,-}w^\sharp_1(x+\xi/\lambda,p,\hq)\;.
\end{equation}
where $\xi=D\hw^\sharp_1(x,p,\hq)$.

Following \cite{ALL} we also note that for any $x'$ close to $x$, we have
$$
\hw_1(x',p,q)\leq q.\hq- \hw_1^\sharp(x',p,\hq)= \hw_1(x,p,\hq)+\hw_1^\sharp(x,p,\hq)- \hw_1^\sharp(x',p,\hq)
$$
because $\hw_1(x,p,q)+\hw_1^\sharp(x,p,\hq)=q.\hq$. Since $\hw_1^\sharp$ has a second order Taylor expansion at
$x$, this gives
\begin{equation}\label{D2+hw1}
-(D\hw^\sharp_1(x,p,\hq), D^2\hw^\sharp_1(x,p,\hq))\in D^{2,+}\hw_1(x,p,q)\;.
\end{equation}

In a symmetric way, if $(x,y, p,q)\in E'_\gamma$ and $\hp=\frac{\partial \hw_2(y,p,q)}{\partial p}$, then 
\begin{equation}\label{D2+w2*}
(D\hw^*_2(y,\hp,q), D^2\hw^*_2(y,\hp,q))\in D^{2,+}w^*_2(y+\xi/\lambda,\hp,q)\;,
\end{equation}
where $\xi=D\hw^*_2(y,\hp,q)$, 
\begin{equation}\label{charcp}
p+\frac{1}{\lambda'}\hp\in \partial^-_\hp w_2^*(y+\xi/\lambda, \hp,q)
\end{equation}
and 
\begin{equation}\label{D2-hw2}
-(D\hw^*_2(y,\hp,q), D^2\hw^*_2(y,\hp,q))\in D^{2,-}\hw_2(y,p,q)\;.
\end{equation}

\noindent {\bf Step 5 : conclusion. } From the previous steps, we know that the set $E'_\gamma$ defined in step 3 has a positive measure for any $\alpha,\; \gamma>0$ sufficiently small.
Hence we can find sequences $\lambda'_n\to+\infty$, $\alpha_n,\; \gamma_n\to 0^+$, $\zeta_n=(\zeta^n_x,\zeta^n_y,\zeta^n_p,\zeta^n_q)\to 0$, $(x_n,y_n, p_n, q_n)$ converging to some 
$(0,0,\barp,\barq)$ such that $(x_n,y_n, p_n, q_n)\in E'_\gamma$ and such that the map
$z_{\zeta_n}$ has a maximum at $(x_n,y_n, p_n, q_n)$.

Let us set 
\begin{equation}\label{Defhpnhqn}
\hp_n=\frac{\partial \hw_2(y_,p_n,q_n)}{\partial p}\;, \qquad \hq_n=\frac{\partial \hw_1(x_,p_n,q_n)}{\partial q}\;,
\end{equation}
$$
(\xi^n_1, X_1^n)=(D\hw^\sharp_1(x_n,p_n,\hq_n), D^2\hw^\sharp_1(x_n,p_n,\hq_n))
$$
and 
$$
(\xi^n_2, X_2^n)=(D\hw^*_2(y_n,\hp_n,q_n), D^2\hw^*_2(y_n,\hp_n,q_n))\;.
$$ 
>From (\ref{charcq}) and (\ref{charcp}) we have 
\begin{equation}\label{hpnhqn}
p_n+\frac{1}{\lambda'_n}\hp_n\in \partial^-_\hp w_2^*(y_n+\xi_2^n/\lambda, \hp_n,q_n)
\;{\rm and}\; q_n+\frac{1}{\lambda'_n}\hq_n\in \partial^+_\hq w_1^\sharp(x_n+\xi^1_n/\lambda, p_n,\hq_n)\;.
\end{equation}
>From (\ref{D2+hw1}) and (\ref{D2-hw2}) we have
$$
-(\xi^n_1, X_1^n)\in D^{2,+}\hw_1(x_n,p_n,q_n)\; {\rm and }\; -(\xi^n_2, X_2^n)\in D^{2,-}\hw_2(x_n,p_n,q_n)\;,
$$
Since furthermore
$(x,y)\to z_{\zeta_n}(x,y,p_n,q_n)$ has a maximum at $(x_n,y_n, p_n, q_n)$, the first and second order optimality conditions imply that
\begin{equation}\label{ishii0n}
(-\xi^n_1, \xi^n_2)=(A+\epsilon A^2)(x_n,y_n)+2\alpha_n(x_n,y_n)+(\zeta^n_x,\zeta^n_y)
\end{equation}
and
\begin{equation}\label{ishii2n}
\left(\frac{1}{\epsilon}+\|A\|\right) I \leq \left(\begin{array}{ll} -X^n_1 & 0\\ 0 & X^n_2\end{array}\right) \leq A+\epsilon A^2+ 2 \alpha_n I
\end{equation}
The left-hand side inequality is due to the fact that $\hw_1$ and $\hw_2$
are semiconvex and semiconcave w.r. to $x$ and $y$ respectively with a modulus $\lambda=\frac{1}{\epsilon}+\|A\|$. 
Using (\ref{D2-w1sharp}) and (\ref{D2+w2*}) gives
\begin{equation}\label{ishii1n}
(\xi^n_1, X_1^n)\in D^{2,-}w_2^\sharp(x_n+\xi^n_1/\lambda,p_n,\hq_n) \; {\rm and}\;
(\xi^n_2, X_2^n) \in D^{2,+}w_1^*(y_n+\xi^n_2/\lambda,\hp_n,q_n) 
\end{equation}

We now note that $(X^n_1)$, $(X^n_2)$, $(\hp_n)$ and $(\hq_n)$ are bounded.
For $(X^n_1)$, $(X^n_2)$ this is an obvious consequence of (\ref{ishii2n}).
For $(\hp_n)$ and $(\hq_n)$ this comes from (\ref{Defhpnhqn}), from the Lipschitz continuity assumption of
$w_2$ and $w_1$ with respect to $p$ and $q$ respectively and from the definition of $\hw_1$ and $\hw_2$. 

We now let $n\to +\infty$. From (\ref{ishii0n}), we have $\xi_1^n,\xi_2^n\to 0$. 
We can assume that $(\hp_n,\hq_n)\to (\hp,\hq)$, $X_1^n\to X_1$ and $X^n_2\to X_2$. Then we have from  (\ref{hpnhqn}), (\ref{ishii1n}) and
 (\ref{ishii2n}) that:
$$
\barp\in \partial^-_\hp w_2^*(0, \hp,\barq)
\;{\rm and}\; \barq\in \partial^+_\hq w_1^\sharp(0, \barp,\hq)\;,
$$
$$
(0, X_1)\in \overline{D^{2,-}}w_2^\sharp(0,\barp,\hq) \; {\rm and}\;
(0, X_2) \in \overline{D^{2,+}}w_1^*(0,\hp,\barq) 
$$
and
$$
\left(\frac{1}{\epsilon}+\|A\|\right) I \leq \left(\begin{array}{ll} -X_1 & 0\\ 0 & X_2\end{array}\right) \leq A+\epsilon A^2\;.
$$
\cq


\end{document}